
\def\Resetstrings{%Clears all strings before processing reference listing.
%   The strings (\Astr, etc.) are fields taken from the database.
%   If the string is present, the appropriate test (\Atest, etc) is set
%   equal to \present---thus allowing the macros to test for the presence 
%   or absence of a field.  All reference processing is done in a local
%   group--the control string names will not conflict with uses
%   elsewhere in the document.
    \def\present{ }\let\bgroup={\let\egroup=}%primitive TeX
    \def\Astr{}\def\astr{}\def\Atest{}\def\atest{}%
    \def\Bstr{}\def\bstr{}\def\Btest{}\def\btest{}%
    \def\Cstr{}\def\cstr{}\def\Ctest{}\def\ctest{}%
    \def\Dstr{}\def\dstr{}\def\Dtest{}\def\dtest{}%
    \def\Estr{}\def\estr{}\def\Etest{}\def\etest{}%
    \def\Fstr{}\def\fstr{}\def\Ftest{}\def\ftest{}%
    \def\Gstr{}\def\gstr{}\def\Gtest{}\def\gtest{}%
    \def\Hstr{}\def\hstr{}\def\Htest{}\def\htest{}%
    \def\Istr{}\def\istr{}\def\Itest{}\def\itest{}%
    \def\Jstr{}\def\jstr{}\def\Jtest{}\def\jtest{}%
    \def\Kstr{}\def\kstr{}\def\Ktest{}\def\ktest{}%
    \def\Lstr{}\def\lstr{}\def\Ltest{}\def\ltest{}%
    \def\Mstr{}\def\mstr{}\def\Mtest{}\def\mtest{}%
    \def\Nstr{}\def\nstr{}\def\Ntest{}\def\ntest{}%
    \def\Ostr{}\def\ostr{}\def\Otest{}\def\otest{}%
    \def\Pstr{}\def\pstr{}\def\Ptest{}\def\ptest{}%
    \def\Qstr{}\def\qstr{}\def\Qtest{}\def\qtest{}%
    \def\Rstr{}\def\rstr{}\def\Rtest{}\def\rtest{}%
    \def\Sstr{}\def\sstr{}\def\Stest{}\def\stest{}%
    \def\Tstr{}\def\tstr{}\def\Ttest{}\def\ttest{}%
    \def\Ustr{}\def\ustr{}\def\Utest{}\def\utest{}%
    \def\Vstr{}\def\vstr{}\def\Vtest{}\def\vtest{}%
    \def\Wstr{}\def\wstr{}\def\Wtest{}\def\wtest{}%
    \def\Xstr{}\def\xstr{}\def\Xtest{}\def\xtest{}%
    \def\Ystr{}\def\ystr{}\def\Ytest{}\def\ytest{}%
}
\Resetstrings

\def\Refformat{%Determines the kind of reference by the presence or
%   absence of certain fields in the database listing, and calls the
%   appropriate macro.
         \if\Jtest\present
             {\if\Vtest\present\journalarticleformat
                  \else\conferencereportformat\fi}
            \else\if\Btest\present\bookarticleformat
               \else\if\Rtest\present\technicalreportformat
                  \else\if\Itest\present\bookformat
                     \else\otherformat\fi\fi\fi\fi}

\def\Rpunct{%Default punctuation control strings if the punctuation
%   is to appear after the citation.  (tib looks for punctuation to
%   precede the incomplete citation in the input document; the TeX
%   output puts it to the left or right depending on the style of citation.)
   \def\Lspace{ }%
   \def\Lperiod{ }%  .
   \def\Lcomma{ }%    ,
   \def\Lquest{ }%     ?
   \def\Lcolon{ }%   :
   \def\Lscolon{ }%   ;
   \def\Lbang{ }%      !
   \def\Lquote{ }%   '
   \def\Lqquote{ }%   "
   \def\Lrquote{ }%    `
   \def\Rspace{}%
   \def\Rperiod{.}%  .
   \def\Rcomma{,}%    ,
   \def\Rquest{?}%     ?
   \def\Rcolon{:}%   :
   \def\Rscolon{;}%   ;
   \def\Rbang{!}%      !
   \def\Rquote{'}%   '
   \def\Rqquote{"}%   "
   \def\Rrquote{`}%    `
   }

\def\Lpunct{%Default punctuation control strings if the punctuation
%   is to appear before the citation.  (tib looks for punctuation to
%   precede the incomplete citation in the input document; the TeX
%   output puts it to the left or right depending on the style of citation.)
   \def\Lspace{}%
   \def\Lperiod{\unskip.}%  .
   \def\Lcomma{\unskip,}%    ,
   \def\Lquest{\unskip?}%     ?
   \def\Lcolon{\unskip:}%   :
   \def\Lscolon{\unskip;}%   ;
   \def\Lbang{\unskip!}%      !
   \def\Lquote{\unskip'}%   '
   \def\Lqquote{\unskip"}%   "
   \def\Lrquote{\unskip`}%    `
   \def\Rspace{\spacefactor=1000}%
   \def\Rperiod{\spacefactor=3000}%  .
   \def\Rcomma{\spacefactor=1250}%    ,
   \def\Rquest{\spacefactor=3000}%     ?
   \def\Rcolon{\spacefactor=2000}%   :
   \def\Rscolon{\spacefactor=1250}%   ;
   \def\Rbang{\spacefactor=3000}%      !
   \def\Rquote{\spacefactor=1000}%   '
   \def\Rqquote{\spacefactor=1000}%   "
   \def\Rrquote{\spacefactor=1000}%    `
   }

\def\Refstd{%Standard control strings for formatting bibliography listings.
     \def\Acomma{\unskip, }%between multiple author names
     \def\Aand{\unskip\ and }%between two author names
     \def\Aandd{\unskip\ and }%between last two of multiple author names
     \def\Ecomma{\unskip, }%between multiple editor names
     \def\Eand{\unskip\ and }%between two editor names
     \def\Eandd{\unskip\ and }%between last two of multiple author names
     \def\acomma{\unskip, }%same for authors of reviewed material
     \def\aand{\unskip\ and }%same for authors of reviewed material
     \def\aandd{\unskip\ and }%same for authors of reviewed material
     \def\ecomma{\unskip, }%same for translators
     \def\eand{\unskip\ and }%same for translators
     \def\eandd{\unskip\ and }%same for translators
     \def\Namecomma{\unskip, }%same for citations using authors' names
     \def\Nameand{\unskip\ and }%same for citations using authors' names
     \def\Nameandd{\unskip\ and }%same for citations using authors' names
     \def\Revcomma{\unskip, }%between last and first name of reversed name
     \def\Initper{.\ }%punctuation after initial
     \def\Initgap{\dimen0=\spaceskip\divide\dimen0 by 2\hskip-\dimen0}%
                           %gap between initials of abbreviated first name
   }

\def\Smallcapsaand{%Smallcaps redefinition of \Aand and \Aandd for \Refstd
     \def\Aand{\unskip\bgroup{\Smallcapsfont\ AND }\egroup}%
     \def\Aandd{\unskip\bgroup{\Smallcapsfont\ AND }\egroup}%
     \def\eand{\unskip\bgroup\Smallcapsfont\ AND \egroup}%
     \def\eandd{\unskip\bgroup\Smallcapsfont\ AND \egroup}%
   }

\def\Smallcapseand{%Smallcaps redefinition of \Eand, \Eeand, etc for Refstd
     \def\Eand{\unskip\bgroup\Smallcapsfont\ AND \egroup}%
     \def\Eandd{\unskip\bgroup\Smallcapsfont\ AND \egroup}%
     \def\aand{\unskip\bgroup\Smallcapsfont\ AND \egroup}%
     \def\aandd{\unskip\bgroup\Smallcapsfont\ AND \egroup}%
   }

%default fonts
   \def\Citefont{}%citations
   \def\ACitefont{}%alternate citations
   \def\Authfont{}%authors
   \def\Titlefont{}%titles
   \def\Tomefont{\sl}%journals or books
   \def\Volfont{}%volume number of journal
   \def\Flagfont{}%citation flag
   \def\Reffont{\rm}%set at beginning of reference listing
   \def\Smallcapsfont{\sevenrm}%small caps
   \def\Flagstyle#1{\hangindent\parindent\indent\hbox to0pt%flag style
       {\hss[{\Flagfont#1}]\kern.5em}\ignorespaces}%        for references

\def\Underlinemark{\vrule height .7pt depth 0pt width 3pc}%for replacing
%   successive listings of identical author(s) by underline (U option
%      in tib format file or -u flag on call).

\def\Citebrackets{\Rpunct%defaults for putting citations in brackets [].
   \def\Lcitemark{\def\Cfont{\Citefont}[\bgroup\Cfont}%mark at left of citation
   \def\Rcitemark{\egroup]}%mark at right of citation
   \def\LAcitemark{\def\Cfont{\ACitefont}\bgroup\ACitefont}%
                         %mark at left of alternate citation
   \def\RAcitemark{\egroup}%mark at right of alternate citation
   \def\LIcitemark{\egroup}%mark at left of insertion in citation
   \def\RIcitemark{\bgroup\Cfont}%mark at right of insertion in citation
   \def\Citehyphen{\egroup--\bgroup\Cfont}%separater for string of citations
   \def\Citecomma{\egroup,\hskip0pt\bgroup\Cfont}%
                                          %separater for multiple citations
   \def\Citebreak{}%mark between parts of citation (e.g. author\Citebreak date)
   }

\def\Citeparen{\Rpunct%defaults for putting citations in parenthesis ().
   \def\Lcitemark{\def\Cfont{\Citefont}(\bgroup\Cfont}%mark at left of citation
   \def\Rcitemark{\egroup)}%mark at right of citation
   \def\LAcitemark{\def\Cfont{\ACitefont}\bgroup\ACitefont}%
                         %mark at left of alternate citation
   \def\RAcitemark{\egroup}%mark at right of alternate citation
   \def\LIcitemark{\egroup}%mark at left of insertion in citation
   \def\RIcitemark{\bgroup\Cfont}%mark at right of insertion in citation
   \def\Citehyphen{\egroup--\bgroup\Cfont}%separater for string of citations
   \def\Citecomma{\egroup,\hskip0pt\bgroup\Cfont}%
                                          %separater for multiple citations
   \def\Citebreak{}%mark between parts of citation (e.g. author\Citebreak date)
   }

\def\Citesuper{\Lpunct%defaults for making superscript citations
   \def\Lcitemark{\def\Cfont{\Citefont}\raise1ex\hbox\bgroup\bgroup\Cfont}%
                         %mark at left of citation
   \def\Rcitemark{\egroup\egroup}%mark at right of citation
   \def\LAcitemark{\def\Cfont{\ACitefont}\bgroup\ACitefont}%
                         %mark at left of alternate citation
   \def\RAcitemark{\egroup}%mark at right of alternate citation
   \def\LIcitemark{\egroup\egroup}%mark at left of insertion in citation
   \def\RIcitemark{\raise1ex\hbox\bgroup\bgroup\Cfont}%
                         %mark at right of insertion in citation
   \def\Citehyphen{\egroup--\bgroup\Cfont}%separater for string of citations
   \def\Citecomma{\egroup,\hskip0pt\bgroup%
      \Cfont}%separater for multiple citations
   \def\Citebreak{}%mark between parts of citation (e.g. author\Citebreak date)
   } 

\def\Citenamedate{\Rpunct%defaults for making name-date citations
   \def\Lcitemark{%mark at left of citation--also sets internal punctuation
      \def\Citebreak{\egroup\ [\bgroup\Citefont}%separater in citation
      \def\Citecomma{\egroup]; %between multiple citations
         \bgroup\let\uchyph=1\Citefont}(\bgroup\let\uchyph=1\Citefont}%
   \def\Rcitemark{\egroup])}%mark at right of citation
   \def\LAcitemark{%mark at left of alternate citation
      \def\Citebreak{\egroup\ [\bgroup\Citefont}\def\Citecomma{\egroup], %
         \bgroup\ACitefont }\bgroup\let\uchyph=1\ACitefont}%
   \def\RAcitemark{\egroup]}%mark at right of alternate citation
  \def\Citehyphen{\egroup--\bgroup\Citefont}%separater for string of citations
   \def\LIcitemark{\egroup}%mark at left of insertion in citation
   \def\RIcitemark{\bgroup\Citefont}%mark at right of insertion in citation
   }
%ams numeric or alpha flag style
%flag, author, title, etc., volume (date) pages, gov't no., other

\Refstd\Citebrackets%set general formats for reference list and citations
\def\Citefont{\bf}\def\Titlefont{\sl}\def\Volfont{\bf}\def\Tomefont{\Reffont}%redefine some fonts

\def\journalarticleformat{\Reffont\let\uchyph=1\parindent=1.25pc\def\Comma{}%
                \sfcode`\.=1000\sfcode`\?=1000\sfcode`\!=1000\sfcode`\:=1000\sfcode`\;=1000\sfcode`\,=1000%\frenchspacing
                \par\vfil\penalty-200\vfilneg%\filbreak
      \if\Ftest\present\Flagstyle\Fstr\fi%
       \if\Atest\present\bgroup\Authfont\Astr\egroup\def\Comma{\unskip, }\fi%
        \if\Ttest\present\Comma\bgroup\Titlefont\Tstr\egroup\def\Comma{, }\fi%
         \if\etest\present\hskip.2em(\bgroup\estr\egroup)\def\Comma{\unskip, }\fi%
          \if\Jtest\present\Comma\bgroup\Tomefont\Jstr\/\egroup\def\Comma{, }\fi%
           \if\Vtest\present\if\Jtest\present\hskip.2em\else\Comma\fi\bgroup\Volfont\Vstr\egroup\def\Comma{, }\fi%
            \if\Dtest\present\hskip.2em(\bgroup\Dstr\egroup)\def\Comma{, }\fi%
             \if\Ptest\present\bgroup, \Pstr\egroup\def\Comma{, }\fi%
              \if\ttest\present\Comma\bgroup\Titlefont\tstr\egroup\def\Comma{, }\fi%
               \if\jtest\present\Comma\bgroup\Tomefont\jstr\/\egroup\def\Comma{, }\fi%
                \if\vtest,\present\if\jtest\present\hskip.2em\else\Comma\fi\bgroup\Volfont\vstr\egroup\def\Comma{, }\fi%
                 \if\dtest\present\hskip.2em(\bgroup\dstr\egroup)\def\Comma{, }\fi%
                  \if\ptest\present\bgroup, \pstr\egroup\def\Comma{, }\fi%
                   \if\Gtest\present{\Comma Gov't ordering no. }\bgroup\Gstr\egroup\def\Comma{, }\fi%
                    \if\Mtest\present\Comma MR \#\bgroup\Mstr\egroup\def\Comma{, }\fi%
                     \if\Otest\present{\Comma\bgroup\Ostr\egroup.}\else{.}\fi%
                      \vskip3ptplus1ptminus1pt}%\smallskip

\def\conferencereportformat{\Reffont\let\uchyph=1\parindent=1.25pc\def\Comma{}%
                \sfcode`\.=1000\sfcode`\?=1000\sfcode`\!=1000\sfcode`\:=1000\sfcode`\;=1000\sfcode`\,=1000%\frenchspacing
                \par\vfil\penalty-200\vfilneg%\filbreak
      \if\Ftest\present\Flagstyle\Fstr\fi%
       \if\Atest\present\bgroup\Authfont\Astr\egroup\def\Comma{\unskip, }\fi%
        \if\Ttest\present\Comma\bgroup\Titlefont\Tstr\egroup\def\Comma{, }\fi%
         \if\Jtest\present\Comma\bgroup\Tomefont\Jstr\/\egroup\def\Comma{, }\fi%
          \if\Ctest\present\Comma\bgroup\Cstr\egroup\def\Comma{, }\fi%
           \if\Dtest\present\hskip.2em(\bgroup\Dstr\egroup)\def\Comma{, }\fi%
            \if\Mtest\present\Comma MR \#\bgroup\Mstr\egroup\def\Comma{, }\fi%
             \if\Otest\present{\Comma\bgroup\Ostr\egroup.}\else{.}\fi%
              \vskip3ptplus1ptminus1pt}%\smallskip

\def\bookarticleformat{\Reffont\let\uchyph=1\parindent=1.25pc\def\Comma{}%
                \sfcode`\.=1000\sfcode`\?=1000\sfcode`\!=1000\sfcode`\:=1000\sfcode`\;=1000\sfcode`\,=1000%\frenchspacing
                \par\vfil\penalty-200\vfilneg%\filbreak
      \if\Ftest\present\Flagstyle\Fstr\fi%
       \if\Atest\present\bgroup\Authfont\Astr\egroup\def\Comma{\unskip, }\fi%
        \if\Ttest\present\Comma\bgroup\Titlefont\Tstr\egroup\def\Comma{, }\fi%
         \if\etest\present\hskip.2em(\bgroup\estr\egroup)\def\Comma{\unskip, }\fi%
          \if\Btest\present\Comma in \bgroup\Tomefont\Bstr\/\egroup\def\Comma{\unskip, }\fi%
           \if\otest\present\ \bgroup\ostr\egroup\def\Comma{, }\fi%
            \if\Etest\present\Comma\bgroup\Estr\egroup\unskip, \ifnum\Ecnt>1eds.\else ed.\fi\def\Comma{, }\fi%
             \if\Stest\present\Comma\bgroup\Sstr\egroup\def\Comma{, }\fi%
              \if\Vtest\present\Comma vol. \bgroup\Vstr\egroup\def\Comma{, }\fi%
               \if\Ntest\present\Comma no. \bgroup\Nstr\egroup\def\Comma{, }\fi%
                \if\Itest\present\Comma\bgroup\Istr\egroup\def\Comma{, }\fi%
                 \if\Ctest\present\Comma\bgroup\Cstr\egroup\def\Comma{, }\fi%
                  \if\Dtest\present\Comma\bgroup\Dstr\egroup\def\Comma{, }\fi%
                   \if\Ptest\present\Comma\Pstr\def\Comma{, }\fi%
                    \if\ttest\present\Comma\bgroup\Titlefont\Tstr\egroup\def\Comma{, }\fi%
                     \if\btest\present\Comma in \bgroup\Tomefont\bstr\egroup\def\Comma{, }\fi%
                       \if\atest\present\Comma\bgroup\astr\egroup\unskip, \if\acnt\present eds.\else ed.\fi\def\Comma{, }\fi%
                        \if\stest\present\Comma\bgroup\sstr\egroup\def\Comma{, }\fi%
                         \if\vtest\present\Comma vol. \bgroup\vstr\egroup\def\Comma{, }\fi%
                          \if\ntest\present\Comma no. \bgroup\nstr\egroup\def\Comma{, }\fi%
                           \if\itest\present\Comma\bgroup\istr\egroup\def\Comma{, }\fi%
                            \if\ctest\present\Comma\bgroup\cstr\egroup\def\Comma{, }\fi%
                             \if\dtest\present\Comma\bgroup\dstr\egroup\def\Comma{, }\fi%
                              \if\ptest\present\Comma\pstr\def\Comma{, }\fi%
                               \if\Gtest\present{\Comma Gov't ordering no. }\bgroup\Gstr\egroup\def\Comma{, }\fi%
                                \if\Mtest\present\Comma MR \#\bgroup\Mstr\egroup\def\Comma{, }\fi%
                                 \if\Otest\present{\Comma\bgroup\Ostr\egroup.}\else{.}\fi%
                                  \vskip3ptplus1ptminus1pt}%\smallskip

\def\bookformat{\Reffont\let\uchyph=1\parindent=1.25pc\def\Comma{}%
                \sfcode`\.=1000\sfcode`\?=1000\sfcode`\!=1000\sfcode`\:=1000\sfcode`\;=1000\sfcode`\,=1000%\frenchspacing
                \par\vfil\penalty-200\vfilneg%\filbreak
      \if\Ftest\present\Flagstyle\Fstr\fi%
       \if\Atest\present\bgroup\Authfont\Astr\egroup\def\Comma{\unskip, }%
            \else\if\Etest\present\bgroup\def\Eand{\Aand}\def\Eandd{\Aandd}\Authfont\Estr\egroup\unskip, \ifnum\Ecnt>1eds.\else ed.\fi\def\Comma{, }%
                  \else\if\Itest\present\bgroup\Authfont\Istr\egroup\def\Comma{, }\fi\fi\fi%
          \if\Ttest\present\Comma\bgroup\Titlefont\Tstr\/\egroup\def\Comma{\unskip, }%
                \else\if\Btest\present\Comma\bgroup\Titlefont\Bstr\/\egroup\def\Comma{\unskip, }\fi\fi%
            \if\otest\present\ \bgroup\ostr\egroup\def\Comma{, }\fi%
             \if\etest\present\hskip.2em(\bgroup\estr\egroup)\def\Comma{\unskip, }\fi%
              \if\Stest\present\Comma\bgroup\Sstr\egroup\def\Comma{, }\fi%
               \if\Vtest\present\Comma vol. \bgroup\Vstr\egroup\def\Comma{, }\fi%
                \if\Ntest\present\Comma no. \bgroup\Nstr\egroup\def\Comma{, }\fi%
                 \if\Atest\present\if\Itest\present
                         \Comma\bgroup\Istr\egroup\def\Comma{\unskip, }\fi%
                      \else\if\Etest\present\if\Itest\present
                              \Comma\bgroup\Istr\egroup\def\Comma{\unskip, }\fi\fi\fi%
                     \if\Ctest\present\Comma\bgroup\Cstr\egroup\def\Comma{, }\fi%
                      \if\Dtest\present\Comma\bgroup\Dstr\egroup\def\Comma{, }\fi%
                       \if\ttest\present\Comma\bgroup\Titlefont\tstr\egroup\def\Comma{, }%
                             \else\if\btest\present\Comma\bgroup\Titlefont\bstr\egroup\def\Comma{, }\fi\fi%
                          \if\stest\present\Comma\bgroup\sstr\egroup\def\Comma{, }\fi%
                           \if\vtest\present\Comma vol. \bgroup\vstr\egroup\def\Comma{, }\fi%
                            \if\ntest\present\Comma no. \bgroup\nstr\egroup\def\Comma{, }\fi%
                             \if\itest\present\Comma\bgroup\istr\egroup\def\Comma{, }\fi%
                              \if\ctest\present\Comma\bgroup\cstr\egroup\def\Comma{, }\fi%
                               \if\dtest\present\Comma\bgroup\dstr\egroup\def\Comma{, }\fi%
                                \if\Gtest\present{\Comma Gov't ordering no. }\bgroup\Gstr\egroup\def\Comma{, }\fi%
                                 \if\Mtest\present\Comma MR \#\bgroup\Mstr\egroup\def\Comma{, }\fi%
                                  \if\Otest\present{\Comma\bgroup\Ostr\egroup.}\else{.}\fi%
                                   \vskip3ptplus1ptminus1pt}%\smallskip

\def\technicalreportformat{\Reffont\let\uchyph=1\parindent=1.25pc\def\Comma{}%
                \sfcode`\.=1000\sfcode`\?=1000\sfcode`\!=1000\sfcode`\:=1000\sfcode`\;=1000\sfcode`\,=1000%\frenchspacing
                \par\vfil\penalty-200\vfilneg%\filbreak
      \if\Ftest\present\Flagstyle\Fstr\fi%
       \if\Atest\present\bgroup\Authfont\Astr\egroup\def\Comma{\unskip, }%
            \else\if\Etest\present\bgroup\def\Eand{\Aand}\def\Eandd{\Aandd}\Authfont\Estr\egroup\unskip, \ifnum\Ecnt>1eds.\else ed.\fi\def\Comma{, }%
                  \else\if\Itest\present\bgroup\Authfont\Istr\egroup\def\Comma{, }\fi\fi\fi%
          \if\Ttest\present\Comma\bgroup\Titlefont\Tstr\egroup\def\Comma{, }\fi%
           \if\Atest\present\if\Itest\present
                   \Comma\bgroup\Istr\egroup\def\Comma{, }\fi%
                \else\if\Etest\present\if\Itest\present
                        \Comma\bgroup\Istr\egroup\def\Comma{, }\fi\fi\fi%
            \if\Rtest\present\Comma\bgroup\Rstr\egroup\def\Comma{, }\fi%
             \if\Ctest\present\Comma\bgroup\Cstr\egroup\def\Comma{, }\fi%
              \if\Dtest\present\Comma\bgroup\Dstr\egroup\def\Comma{, }\fi%
               \if\ttest\present\Comma\bgroup\Titlefont\tstr\egroup\def\Comma{, }\fi%
                \if\itest\present\Comma\bgroup\istr\egroup\def\Comma{, }\fi%
                 \if\rtest\present\Comma\bgroup\rstr\egroup\def\Comma{, }\fi%
                  \if\ctest\present\Comma\bgroup\cstr\egroup\def\Comma{, }\fi%
                   \if\dtest\present\Comma\bgroup\dstr\egroup\def\Comma{, }\fi%
                    \if\Gtest\present{\Comma Gov't ordering no. }\bgroup\Gstr\egroup\def\Comma{, }\fi%
                     \if\Mtest\present\Comma MR \#\bgroup\Mstr\egroup\def\Comma{, }\fi%
                      \if\Otest\present{\Comma\bgroup\Ostr\egroup.}\else{.}\fi%
                       \vskip3ptplus1ptminus1pt}%\smallskip

\def\otherformat{\Reffont\let\uchyph=1\parindent=1.25pc\def\Comma{}%
                \sfcode`\.=1000\sfcode`\?=1000\sfcode`\!=1000\sfcode`\:=1000\sfcode`\;=1000\sfcode`\,=1000%\frenchspacing
                \par\vfil\penalty-200\vfilneg%\filbreak
      \if\Ftest\present\Flagstyle\Fstr\fi%
       \if\Atest\present\bgroup\Authfont\Astr\egroup\def\Comma{\unskip, }%
            \else\if\Etest\present\bgroup\def\Eand{\Aand}\def\Eandd{\Aandd}\Authfont\Estr\egroup\unskip, \ifnum\Ecnt>1eds.\else ed.\fi\def\Comma{, }%
                  \else\if\Itest\present\bgroup\Authfont\Istr\egroup\def\Comma{, }\fi\fi\fi%
          \if\Ttest\present\Comma\bgroup\Titlefont\Tstr\egroup\def\Comma{, }\fi%
            \if\Atest\present\if\Itest\present
                    \Comma\bgroup\Istr\egroup\def\Comma{, }\fi%
                 \else\if\Etest\present\if\Itest\present
                         \Comma\bgroup\Istr\egroup\def\Comma{, }\fi\fi\fi%
                 \if\Ctest\present\Comma\bgroup\Cstr\egroup\def\Comma{, }\fi%
                  \if\Dtest\present\Comma\bgroup\Dstr\egroup\def\Comma{, }\fi%
                   \if\Gtest\present{\Comma Gov't ordering no. }\bgroup\Gstr\egroup\def\Comma{, }\fi%
                    \if\Mtest\present\Comma MR \#\bgroup\Mstr\egroup\def\Comma{, }\fi%
                     \if\Otest\present{\Comma\bgroup\Ostr\egroup.}\else{.}\fi%
                      \vskip3ptplus1ptminus1pt}%\smallskip
%ams numeric flag style
%flag followed by period

\def\Flagstyle#1{\Flagfont#1. }%flag style
\message {)}\message {DOCUMENT TEXT}
\message {REFERENCE FORMATTING FILES:}
%ams alpha flag style
%flag enclosed in brackets

\def\Flagfont{\bf}\def\Citefont{\bf}\def\ACitefont{\bf}%redefine some fonts
\def\Flagstyle#1{[{\Flagfont#1}] }%flag style
%ams numeric or alpha flag style
%flag, author, title, etc., volume (date) pages, gov't no., other

\Refstd\Citebrackets%set general formats for reference list and citations
\def\Citefont{\bf}\def\Titlefont{\sl}\def\Volfont{\bf}\def\Tomefont{\Reffont}%redefine some fonts

\def\journalarticleformat{\Reffont\let\uchyph=1\parindent=1.25pc\def\Comma{}%
                \sfcode`\.=1000\sfcode`\?=1000\sfcode`\!=1000\sfcode`\:=1000\sfcode`\;=1000\sfcode`\,=1000%\frenchspacing
                \par\vfil\penalty-200\vfilneg%\filbreak
      \if\Ftest\present\Flagstyle\Fstr\fi%
       \if\Atest\present\bgroup\Authfont\Astr\egroup\def\Comma{\unskip, }\fi%
        \if\Ttest\present\Comma\bgroup\Titlefont\Tstr\egroup\def\Comma{, }\fi%
         \if\etest\present\hskip.2em(\bgroup\estr\egroup)\def\Comma{\unskip, }\fi%
          \if\Jtest\present\Comma\bgroup\Tomefont\Jstr\/\egroup\def\Comma{, }\fi%
           \if\Vtest\present\if\Jtest\present\hskip.2em\else\Comma\fi\bgroup\Volfont\Vstr\egroup\def\Comma{, }\fi%
            \if\Dtest\present\hskip.2em(\bgroup\Dstr\egroup)\def\Comma{, }\fi%
             \if\Ptest\present\bgroup, \Pstr\egroup\def\Comma{, }\fi%
              \if\ttest\present\Comma\bgroup\Titlefont\tstr\egroup\def\Comma{, }\fi%
               \if\jtest\present\Comma\bgroup\Tomefont\jstr\/\egroup\def\Comma{, }\fi%
                \if\vtest,\present\if\jtest\present\hskip.2em\else\Comma\fi\bgroup\Volfont\vstr\egroup\def\Comma{, }\fi%
                 \if\dtest\present\hskip.2em(\bgroup\dstr\egroup)\def\Comma{, }\fi%
                  \if\ptest\present\bgroup, \pstr\egroup\def\Comma{, }\fi%
                   \if\Gtest\present{\Comma Gov't ordering no. }\bgroup\Gstr\egroup\def\Comma{, }\fi%
                    \if\Mtest\present\Comma MR \#\bgroup\Mstr\egroup\def\Comma{, }\fi%
                     \if\Otest\present{\Comma\bgroup\Ostr\egroup.}\else{.}\fi%
                      \vskip3ptplus1ptminus1pt}%\smallskip

\def\conferencereportformat{\Reffont\let\uchyph=1\parindent=1.25pc\def\Comma{}%
                \sfcode`\.=1000\sfcode`\?=1000\sfcode`\!=1000\sfcode`\:=1000\sfcode`\;=1000\sfcode`\,=1000%\frenchspacing
                \par\vfil\penalty-200\vfilneg%\filbreak
      \if\Ftest\present\Flagstyle\Fstr\fi%
       \if\Atest\present\bgroup\Authfont\Astr\egroup\def\Comma{\unskip, }\fi%
        \if\Ttest\present\Comma\bgroup\Titlefont\Tstr\egroup\def\Comma{, }\fi%
         \if\Jtest\present\Comma\bgroup\Tomefont\Jstr\/\egroup\def\Comma{, }\fi%
          \if\Ctest\present\Comma\bgroup\Cstr\egroup\def\Comma{, }\fi%
           \if\Dtest\present\hskip.2em(\bgroup\Dstr\egroup)\def\Comma{, }\fi%
            \if\Mtest\present\Comma MR \#\bgroup\Mstr\egroup\def\Comma{, }\fi%
             \if\Otest\present{\Comma\bgroup\Ostr\egroup.}\else{.}\fi%
              \vskip3ptplus1ptminus1pt}%\smallskip

\def\bookarticleformat{\Reffont\let\uchyph=1\parindent=1.25pc\def\Comma{}%
                \sfcode`\.=1000\sfcode`\?=1000\sfcode`\!=1000\sfcode`\:=1000\sfcode`\;=1000\sfcode`\,=1000%\frenchspacing
                \par\vfil\penalty-200\vfilneg%\filbreak
      \if\Ftest\present\Flagstyle\Fstr\fi%
       \if\Atest\present\bgroup\Authfont\Astr\egroup\def\Comma{\unskip, }\fi%
        \if\Ttest\present\Comma\bgroup\Titlefont\Tstr\egroup\def\Comma{, }\fi%
         \if\etest\present\hskip.2em(\bgroup\estr\egroup)\def\Comma{\unskip, }\fi%
          \if\Btest\present\Comma in \bgroup\Tomefont\Bstr\/\egroup\def\Comma{\unskip, }\fi%
           \if\otest\present\ \bgroup\ostr\egroup\def\Comma{, }\fi%
            \if\Etest\present\Comma\bgroup\Estr\egroup\unskip, \ifnum\Ecnt>1eds.\else ed.\fi\def\Comma{, }\fi%
             \if\Stest\present\Comma\bgroup\Sstr\egroup\def\Comma{, }\fi%
              \if\Vtest\present\Comma vol. \bgroup\Vstr\egroup\def\Comma{, }\fi%
               \if\Ntest\present\Comma no. \bgroup\Nstr\egroup\def\Comma{, }\fi%
                \if\Itest\present\Comma\bgroup\Istr\egroup\def\Comma{, }\fi%
                 \if\Ctest\present\Comma\bgroup\Cstr\egroup\def\Comma{, }\fi%
                  \if\Dtest\present\Comma\bgroup\Dstr\egroup\def\Comma{, }\fi%
                   \if\Ptest\present\Comma\Pstr\def\Comma{, }\fi%
                    \if\ttest\present\Comma\bgroup\Titlefont\Tstr\egroup\def\Comma{, }\fi%
                     \if\btest\present\Comma in \bgroup\Tomefont\bstr\egroup\def\Comma{, }\fi%
                       \if\atest\present\Comma\bgroup\astr\egroup\unskip, \if\acnt\present eds.\else ed.\fi\def\Comma{, }\fi%
                        \if\stest\present\Comma\bgroup\sstr\egroup\def\Comma{, }\fi%
                         \if\vtest\present\Comma vol. \bgroup\vstr\egroup\def\Comma{, }\fi%
                          \if\ntest\present\Comma no. \bgroup\nstr\egroup\def\Comma{, }\fi%
                           \if\itest\present\Comma\bgroup\istr\egroup\def\Comma{, }\fi%
                            \if\ctest\present\Comma\bgroup\cstr\egroup\def\Comma{, }\fi%
                             \if\dtest\present\Comma\bgroup\dstr\egroup\def\Comma{, }\fi%
                              \if\ptest\present\Comma\pstr\def\Comma{, }\fi%
                               \if\Gtest\present{\Comma Gov't ordering no. }\bgroup\Gstr\egroup\def\Comma{, }\fi%
                                \if\Mtest\present\Comma MR \#\bgroup\Mstr\egroup\def\Comma{, }\fi%
                                 \if\Otest\present{\Comma\bgroup\Ostr\egroup.}\else{.}\fi%
                                  \vskip3ptplus1ptminus1pt}%\smallskip

\def\bookformat{\Reffont\let\uchyph=1\parindent=1.25pc\def\Comma{}%
                \sfcode`\.=1000\sfcode`\?=1000\sfcode`\!=1000\sfcode`\:=1000\sfcode`\;=1000\sfcode`\,=1000%\frenchspacing
                \par\vfil\penalty-200\vfilneg%\filbreak
      \if\Ftest\present\Flagstyle\Fstr\fi%
       \if\Atest\present\bgroup\Authfont\Astr\egroup\def\Comma{\unskip, }%
            \else\if\Etest\present\bgroup\def\Eand{\Aand}\def\Eandd{\Aandd}\Authfont\Estr\egroup\unskip, \ifnum\Ecnt>1eds.\else ed.\fi\def\Comma{, }%
                  \else\if\Itest\present\bgroup\Authfont\Istr\egroup\def\Comma{, }\fi\fi\fi%
          \if\Ttest\present\Comma\bgroup\Titlefont\Tstr\/\egroup\def\Comma{\unskip, }%
                \else\if\Btest\present\Comma\bgroup\Titlefont\Bstr\/\egroup\def\Comma{\unskip, }\fi\fi%
            \if\otest\present\ \bgroup\ostr\egroup\def\Comma{, }\fi%
             \if\etest\present\hskip.2em(\bgroup\estr\egroup)\def\Comma{\unskip, }\fi%
              \if\Stest\present\Comma\bgroup\Sstr\egroup\def\Comma{, }\fi%
               \if\Vtest\present\Comma vol. \bgroup\Vstr\egroup\def\Comma{, }\fi%
                \if\Ntest\present\Comma no. \bgroup\Nstr\egroup\def\Comma{, }\fi%
                 \if\Atest\present\if\Itest\present
                         \Comma\bgroup\Istr\egroup\def\Comma{\unskip, }\fi%
                      \else\if\Etest\present\if\Itest\present
                              \Comma\bgroup\Istr\egroup\def\Comma{\unskip, }\fi\fi\fi%
                     \if\Ctest\present\Comma\bgroup\Cstr\egroup\def\Comma{, }\fi%
                      \if\Dtest\present\Comma\bgroup\Dstr\egroup\def\Comma{, }\fi%
                       \if\ttest\present\Comma\bgroup\Titlefont\tstr\egroup\def\Comma{, }%
                             \else\if\btest\present\Comma\bgroup\Titlefont\bstr\egroup\def\Comma{, }\fi\fi%
                          \if\stest\present\Comma\bgroup\sstr\egroup\def\Comma{, }\fi%
                           \if\vtest\present\Comma vol. \bgroup\vstr\egroup\def\Comma{, }\fi%
                            \if\ntest\present\Comma no. \bgroup\nstr\egroup\def\Comma{, }\fi%
                             \if\itest\present\Comma\bgroup\istr\egroup\def\Comma{, }\fi%
                              \if\ctest\present\Comma\bgroup\cstr\egroup\def\Comma{, }\fi%
                               \if\dtest\present\Comma\bgroup\dstr\egroup\def\Comma{, }\fi%
                                \if\Gtest\present{\Comma Gov't ordering no. }\bgroup\Gstr\egroup\def\Comma{, }\fi%
                                 \if\Mtest\present\Comma MR \#\bgroup\Mstr\egroup\def\Comma{, }\fi%
                                  \if\Otest\present{\Comma\bgroup\Ostr\egroup.}\else{.}\fi%
                                   \vskip3ptplus1ptminus1pt}%\smallskip

\def\technicalreportformat{\Reffont\let\uchyph=1\parindent=1.25pc\def\Comma{}%
                \sfcode`\.=1000\sfcode`\?=1000\sfcode`\!=1000\sfcode`\:=1000\sfcode`\;=1000\sfcode`\,=1000%\frenchspacing
                \par\vfil\penalty-200\vfilneg%\filbreak
      \if\Ftest\present\Flagstyle\Fstr\fi%
       \if\Atest\present\bgroup\Authfont\Astr\egroup\def\Comma{\unskip, }%
            \else\if\Etest\present\bgroup\def\Eand{\Aand}\def\Eandd{\Aandd}\Authfont\Estr\egroup\unskip, \ifnum\Ecnt>1eds.\else ed.\fi\def\Comma{, }%
                  \else\if\Itest\present\bgroup\Authfont\Istr\egroup\def\Comma{, }\fi\fi\fi%
          \if\Ttest\present\Comma\bgroup\Titlefont\Tstr\egroup\def\Comma{, }\fi%
           \if\Atest\present\if\Itest\present
                   \Comma\bgroup\Istr\egroup\def\Comma{, }\fi%
                \else\if\Etest\present\if\Itest\present
                        \Comma\bgroup\Istr\egroup\def\Comma{, }\fi\fi\fi%
            \if\Rtest\present\Comma\bgroup\Rstr\egroup\def\Comma{, }\fi%
             \if\Ctest\present\Comma\bgroup\Cstr\egroup\def\Comma{, }\fi%
              \if\Dtest\present\Comma\bgroup\Dstr\egroup\def\Comma{, }\fi%
               \if\ttest\present\Comma\bgroup\Titlefont\tstr\egroup\def\Comma{, }\fi%
                \if\itest\present\Comma\bgroup\istr\egroup\def\Comma{, }\fi%
                 \if\rtest\present\Comma\bgroup\rstr\egroup\def\Comma{, }\fi%
                  \if\ctest\present\Comma\bgroup\cstr\egroup\def\Comma{, }\fi%
                   \if\dtest\present\Comma\bgroup\dstr\egroup\def\Comma{, }\fi%
                    \if\Gtest\present{\Comma Gov't ordering no. }\bgroup\Gstr\egroup\def\Comma{, }\fi%
                     \if\Mtest\present\Comma MR \#\bgroup\Mstr\egroup\def\Comma{, }\fi%
                      \if\Otest\present{\Comma\bgroup\Ostr\egroup.}\else{.}\fi%
                       \vskip3ptplus1ptminus1pt}%\smallskip

\def\otherformat{\Reffont\let\uchyph=1\parindent=1.25pc\def\Comma{}%
                \sfcode`\.=1000\sfcode`\?=1000\sfcode`\!=1000\sfcode`\:=1000\sfcode`\;=1000\sfcode`\,=1000%\frenchspacing
                \par\vfil\penalty-200\vfilneg%\filbreak
      \if\Ftest\present\Flagstyle\Fstr\fi%
       \if\Atest\present\bgroup\Authfont\Astr\egroup\def\Comma{\unskip, }%
            \else\if\Etest\present\bgroup\def\Eand{\Aand}\def\Eandd{\Aandd}\Authfont\Estr\egroup\unskip, \ifnum\Ecnt>1eds.\else ed.\fi\def\Comma{, }%
                  \else\if\Itest\present\bgroup\Authfont\Istr\egroup\def\Comma{, }\fi\fi\fi%
          \if\Ttest\present\Comma\bgroup\Titlefont\Tstr\egroup\def\Comma{, }\fi%
            \if\Atest\present\if\Itest\present
                    \Comma\bgroup\Istr\egroup\def\Comma{, }\fi%
                 \else\if\Etest\present\if\Itest\present
                         \Comma\bgroup\Istr\egroup\def\Comma{, }\fi\fi\fi%
                 \if\Ctest\present\Comma\bgroup\Cstr\egroup\def\Comma{, }\fi%
                  \if\Dtest\present\Comma\bgroup\Dstr\egroup\def\Comma{, }\fi%
                   \if\Gtest\present{\Comma Gov't ordering no. }\bgroup\Gstr\egroup\def\Comma{, }\fi%
                    \if\Mtest\present\Comma MR \#\bgroup\Mstr\egroup\def\Comma{, }\fi%
                     \if\Otest\present{\Comma\bgroup\Ostr\egroup.}\else{.}\fi%
                      \vskip3ptplus1ptminus1pt}%\smallskip

\message {)}\message{)}\message {DOCUMENT TEXT}
%Does not contain references tibbed.
\documentstyle{amsppt}
\magnification =\magstep1

\hcorrection{.25in}

\NoBlackBoxes
\def\D{\Delta}
\predefine\Sec{\S}
\redefine\L{\Lambda}
\def\O{\Omega}
\def\P{\Phi}
\redefine\S{\Sigma}
\def\a{\alpha}
\def\b{\beta}
\def\d{\delta}
\def\e{\varepsilon}
\def\g{\gamma}
\def\l{\lambda}
\def\o{\omega}

\def\p{\psi}
\def\s{\sigma}
\def\t{\chi}
\def\T{\theta}

\def\C{\Bbb C}

\def\Z{\Bbb Z}

\def\A{\text{\bf A}}
\def\B{\text{\bf B}}
\def\G{\text{\bf G}}
\def\M{\text{\bf M}}
\def\N{\bold N}
\def\PP{\text{\bf P}}
\def\TT{\text{\bf T}}
\def\U{\text{\bf U}}
\def\aa{\text{\bf a}}
\def\x{\text{\bf x}}

\def\Hom{\operatorname{Hom}}
\def\Ind{\operatorname{Ind}}
\def\ind{\operatorname{ind}}
\def\mod{\operatorname{mod}}
\def\Res{\operatorname{Res}}

\def\cc{{C_c^\infty}}

\def\tot{{\tilde{\omega}_\chi}}
\def\DD{{\Cal D}}
\def\PPP{{\Cal P}}
\loadeurm
\topmatter
\title
On local coefficients for non-generic representations of some classical groups
\endtitle
\author Solomon Friedberg  and David Goldberg\endauthor
%\affil University of California, Santa Cruz \\
%Purdue University \endaffil
\address 
Solomon Friedberg,
Department of Mathematics,
University of California Santa Cruz,
Santa Cruz, CA 95064
\endaddress
\email
friedbe\@cats.ucsc.edu
\endemail
\address  
David Goldberg,
Department of Mathematics,
Purdue University,
West Lafayette, IN 47907
\endaddress
\email
goldberg\@math.purdue.edu
\endemail
\rightheadtext{On local coefficients for non-generic representations}
\abstract This paper is concerned with representations
of split orthogonal and quasi-split unitary groups over a 
nonarchimedean local field which
are not generic, but which support a unique model of a different
kind, the generalized Bessel model.  The properties 
of the Bessel models under induction are studied, and an analogue of
Rodier's theorem concerning the induction of Whittaker models
is proved for Bessel models which are
minimal in a suitable sense.  The holomorphicity in the induction 
parameter of the Bessel functional is established.
Last, local coefficients are defined for each 
irreducible supercuspidal representation
which carries a Bessel functional and also for a certain
component of each representation parabolically induced from 
such a supercuspidal.
\endabstract
%\leftheadtext{Solomon Friedberg and David Goldberg}

\thanks{Friedberg was partially supported by National Security
Administration Grant MDA904-95-H-1053.  Goldberg was partially
supported by National Science Foundation Fellowship DMS-9206246 and
National Science Foundation Career Grant DMS-9501868.  The research of
both authors at MSRI was supported in part by National Science
Foundation Grant DMS-9022140.}
\endthanks

\endtopmatter
\document

\subheading{Introduction}
$L$-functions are a central object of study in representation theory 
and number theory.  Over a global field, one has the Langlands
Conjectures, which assert in particular the meromorphic continuation
and functional equation of a class of Euler products.   Over
a local field one has additional conjectures due to 
Langlands,
expressing the Plancherel measure arithmetically as the ratio
of certain local $L$-functions and root numbers.

In many cases these conjectures have been established by
Shahidi\Lspace \Lcitemark Shab\Citecomma
Shac\Citecomma
Shad\Rcitemark \Rspace{}, 
following a path laid out by Langlands\Lspace \Lcitemark Lana\Rcitemark \Rspace{}.  The framework
for Shahidi's work is the study of Eisenstein series or their
local analogues, induced representations.  One knows the
continuation of these Eisenstein series due to Langlands
\Lcitemark Lanb\Rcitemark \Rspace{}.
Langlands also showed that the constant coefficients 
of the Eisenstein series may be expressed in terms of
local intertwining operators which are almost everywhere
quotients of certain $L$-functions.  It remains to
study these intertwining operators for the finite set
of `bad' places.  If the inducing data is generic,
that is, admits a Whittaker model, then Shahidi has
succeeded in relating them to
local $L$-functions.  Thus the careful study of
the Eisenstein series, both local and global, affords a
proof of certain of the Langlands conjectures for these $L$-functions.

The aim of this work is to suggest that the Langlands-Shahidi
method may be extended beyond the generic spectrum by the
use of other models.  The Whittaker model is unique (an
irreducible admissible representation admits at most one such
model up to scalars).  In this paper we study the 
properties of local representations of split orthogonal groups and 
quasi-split unitary groups which are not generic, but which support a  
unique model of a different kind, the generalized Bessel model.
These models involve a character of a proper subgroup of
the unipotent radical of a Borel subgroup, but transform
under a reductive group of some, in general non-zero,
rank.  The uniqueness of the models has been proved by S. Rallis
\Lcitemark Ral\Rcitemark \Rspace{}
in the orthogonal case, but as the argument has not yet
been written out in full detail in the unitary case
we make it a hypothesis throughout the paper.

We first study the properties of Bessel models under induction,
and prove an analogue of Rodier's Theorem\Lspace \Lcitemark Rodb\Rcitemark \Rspace{} concerning the
induction of Whittaker models.  Our analogue, Theorem 2.1, states that
if one parabolically induces a representation with a
Bessel model of minimal rank, or more generally one which is 
minimal in the sense of Definition 1.5 below, then the induced representation
has a unique Bessel model of the same rank and compatible type.
In the case of rank 0, we recover Rodier's
theorem.  To carry out the proof we use Bruhat's extension\Lspace \Lcitemark Bru\Rcitemark \Rspace{}
of Mackey theory and investigate precisely which double
cosets of the appropriate type may support a functional
with the desired equivariance property.  We show that 
there is a unique such double coset by an extensive
combinatorial argument.

Next, we establish the holomorphicity of the Bessel functional
which arises from one which is minimal by parabolic induction
of the underlying representation.
Our approach is based on Bernstein's theorem\Lspace \Lcitemark Ber\Rcitemark \Rspace{}, which uses
uniqueness to conclude meromorphicity under some regularity
hypotheses, and Banks's extension\Lspace \Lcitemark Ban\Rcitemark \Rspace{}, which allows one to prove
holomorphicity as well.  We show in Theorem 3.6 that there is a
non-zero Bessel
functional $\Lambda(\nu,\pi)$, attached to an irreducible
admissible representation $\pi$ of the Levi subgroup $M$
and a parameter $\nu$ in the complexified dual of the
Lie algebra of the split component of $M$, 
which is holomorphic in $\nu$.

If $\pi$ is supercuspidal and has a Bessel model, 
or more generally if $\pi$ is irreducible and 
carries a Bessel model corresponding to a 
minimal Bessel model 
of the supercuspidal from which it is induced,
these results allow us to establish the existence of
a local coefficient.  In the generic case, such a
local coefficient was crucial 
for Shahidi's study of the intertwining operators
and of the relation between Plancherel measures and $L$--functions;
see Shahidi\Lspace \Lcitemark Shad\Rcitemark \Rspace{}. 
Let $A(\nu,\pi,w)$ denote the standard intertwining
operator attached to inducing data $\nu,\pi$ and Weyl
group element representative $w$ (see (3.4) below).
We shall prove (cf. Theorem 3.8):

\proclaim{Theorem} Let $\pi$ be an irreducible representation
of $M$ which is a component of the representation parabolically
induced from
an irreducible admissible supercuspidal
representation $\rho$ of a parabolic subgroup of $M.$   
Suppose that $\pi$
carries a Bessel model corresponding 
(in the sense of Theorem 2.1) to a minimal Bessel model of $\rho$.
For each $\tilde w$ in the Weyl group, choose a representative
$w$ for $\tilde w.$  Then
there is a complex number $C(\nu,\pi,w)$
so that
$$\L(\nu,\pi)=C(\nu,\pi,w)\L(\tilde w\nu,\tilde w\pi)A(\nu,\pi,w).$$
Moreover, the function $\nu\mapsto C(\nu,\pi,w)$ is meromorphic
and depends only on the class of $\pi$ and the choice of
the representative $w.$
\endproclaim

We call $C(\nu,\pi,w)$ the {\it local coefficient} attached
to $\pi$, $\nu$, and $w$.    We then show in Corollary
3.9  that the local coefficients
behave as expected with respect to the Langlands decomposition
of the intertwining operators.  This 
generalizes a property of the local coefficients
introduced by Shahidi\Lspace \Lcitemark Shaa\Rcitemark \Rspace{} in the generic case.

The authors wish to thank D.\ Ramakrishnan, who suggested we begin
our study of using the Bessel functional to extend the Langlands-Shahidi
method by seeking to generalize Rodier's
theorem, and F.\ Shahidi for his interest in this project. 
The authors also express their appreciation to W.\ Banks,
D.\ Bump,  M.\ Furusawa,
D.\ Ginzburg, S.\ Rallis, F.\ Shahidi, and D.\ Soudry for 
helpful discussions of their work.  Finally, the authors wish to
thank the Mathematical Sciences Research Institute and the organizers
of the 1994-95 Special Year in Automorphic Forms.  The authors
met and began this work while both were in residence there; 
without their support this project would not have started.

\subheading{\Sec 1 Preliminaries on Bessel Models}
In this section we recall the notion of a Bessel model
following \Lspace \Lcitemark Ral\Rcitemark \Rspace{} and
\Lcitemark GPR\Rcitemark \Rspace{}, and review some properties
of such models.
Let $F$ be a nonarchimedean local field of characteristic zero.  Let
$\G$ be one of the classical groups
$SO_{2r+1}, U_{2r+1}, U_{2r},$ or $SO_{2r},$ defined over $F.$
We assume that the orthogonal groups are split, and that the
unitary groups are quasi-split, and split over a
quadratic extension $E/F.$
Let $r_0=2r$ if $\G=U_{2r}$ or $SO_{2r},$
and $r_0=2r+1$ otherwise.
Denote by $\B=\TT\U$ the Borel subgroup of $\G,$ where
$\TT$ contains the maximal split subtorus of diagonal elements,
and $\U$ is the subgroup of upper triangular unipotent matrices in $\G.$
We use  $G$ to denote the $F$\snug-rational points of $\G,$
and use this notational convention for other algebraic
groups defined over $F.$

Denote by $\P(\G,\TT)$ the root system of $\G$
with respect to $\TT.$  We choose the ordering
on the roots corresponding to our choice of Borel subgroup.
Let $W=W(\G,\TT_d)$ be the Weyl
group of $\G$ with respect to the maximal split
subtorus $\TT_d$  of $\TT.$ Thus, $W=N_G(\TT_d)/\TT.$
Then,
$$W\simeq\cases S_r\ltimes\Z_2^r&\text{if }\G\neq SO_{2r}\\
S_r\ltimes \Z_2^{r-1}&\text{if }\G=SO_{2r}.\endcases$$
(See\Lspace \Lcitemark Gola\Citecomma
Golb\Rcitemark \Rspace{} for a more
explicit description of $\TT$ and $W.$)
Here we will denote all elements of $W$
as permutations on $r_0$ letters.
Thus, the permutation $(ij)\in S_r$ corresponds
to the permutation $(ij)(r_0+1-j\, r_0+1-i)$ in $S_{r_0}.$
Similarly, the sign change $c_i$ which generates that $i$\snug-th
copy of $\Z_2$ corresponds to the permutation $(i\,r_0+1-i)$ in
$S_{r_0}.$

Fix an $\ell<r$ and let $\ell_0=r_0-2\ell.$
Let $\U_\ell$ be the subgroup of $\U$
consisting of matrices whose middle $\ell_0\times\ell_0$
block is the identity matrix.
For $1\leq i\leq \ell,$
let $\p_i$ be a non-trivial additive character of $F$
if $\G$ is orthogonal, and let
$\p_i$ be the composition of such
a character with $Tr_{E/F}$
if $\G$ is unitary.
We let $\aa=(a_1,a_2,\dots,a_{\ell_0})\in F^{\ell_0}$
if $\G$ is orthogonal, and let $\aa\in E^{\ell_0}$
if $\G$ is unitary.
Then define $\p_{\ell,a_j}$ by
$\p_{\ell,a_j}(x)=\p_{\ell}(a_jx).$
Now  define a character of $U_{\ell}$
by
$$\t\left(\left(u_{ij}\right)\right)=\prod_{i=1}^{\ell-1}\p_i(u_{i,i+1})\,\prod_{j=1}^{\ell_0}\p_{\ell,a_j}(u_{\ell,\ell+j}).$$
Let $$\M_\ell=\left\{\pmatrix I_\ell\\ &g\\&&I_\ell\endpmatrix\in \G\right\}.$$
Note that $\M_\ell\subset N_{\G}(\U_\ell).$
If $g\in M_\ell,$ then define $\t^g$
by $\t^g(u)=\t(g^{-1}ug).$
We let 
$$M_\t=\left\{g\in M_\ell\,|\, \t^g=\t\right\}.$$
Let $R_\t=M_\t U_\ell.$  Suppose that $\o$ is
an irreducible admissible representation of $M_\t.$
(We will denote this by $\o\in\Cal E(M_\t).$)
Let $\o_\t=\o\otimes\t$ be the
associated representation on $R_\t.$

\definition{Definition 1.1}
We say that two characters $\t_1$ and $\t_2$ of $U_\ell$
defined as above are {\it equivalent} if $\t_1=\t_2^g$
for some $g\in N_G(U_\ell).$
\enddefinition

The following result is a consequence of Witt's Theorem.

\proclaim{Lemma 1.2}  Any character $\t$ of $U_\ell$
which is defined as above, is equivalent to one
for which  $\aa=(\d,0,0\dots,1),$ for some $\d.$
\endproclaim

From now on we assume for convenience that $\t$ is given as in Lemma 1.2.

We let $\ell_1=\left[\frac{\ell_0}{2}\right]=r-\ell.$

\definition{Definition 1.3}  Suppose that $\tau$
is an admissible representation of $G$.
We say that $\tau$ has an $\o_\t$\snug-{\it Bessel model}
(or {\it a Bessel model with respect to $\o_\t$}) if
$\Hom_G(\tau,\Ind_{R_\t}^G(\o_\t))\neq 0.$
If $\t$ is a character of $U_\ell,$ and $\ell_1$ is defined
as above,
then we say that $\tau$ has a {\it rank $\ell_1$}
Bessel model.
\enddefinition

\remark{Remarks}
\roster
\item
By Frobenius reciprocity\Lspace \Lcitemark BeZb\Rcitemark \Rspace{}, we have
$$\Hom_G(\tau,\Ind_{R_\t}^G(\o_\t))\simeq 
\Hom_{M_\t}(\tau_{U_{\ell,\t}},\o_\t),$$
where $\tau_{U_{\ell,\t}}$ is the $\t$\snug-twisted
Jacquet module of $\tau$ with respect to $U_\ell$\Lspace \Lcitemark BeZa\Rcitemark \Rspace{}.
Thus, the non-vanishing of
$\tau_{U_{\ell,\t}},$ for some $\ell$ and $\t,$
would imply that $\tau$ has a rank $\ell_1$ Bessel model
with respect to some $\o_\t.$
\item A Whittaker model\Lspace \Lcitemark Roda\Citecomma
Rodb\Rcitemark \Rspace{}
is a rank zero Bessel model.
\item  One can make these definitions for any choice
of Borel subgroup.  We choose the standard one for
convenience, but we will need to use others in the sequel.
\item
When $\G=SO_{2r+1},$ Rallis has shown that every irreducible
admissible representation of $\G(F)$ has a Bessel model for some choice
of $\t$ and $\o.$  
\item
Suppose that $\tau$ is irreducible, $\o\in\Cal E(M_\t),$
and  $\l:V_{\tau}\rightarrow V_{\o}=V_{\o_\t}$
satisfies $\l(\tau(x)v)=\d_{R_\t}(x)^{1/2}\o_\t(x)\l(v),$ for all
$x\in R_\t$ and $v\in V_\tau.$  (Such a 
$\l$ is called a Bessel functional.) 
Let $v\in V_{\tau}$ and set $B_v(g)=\l(\tau(g)v).$
Then the map $v\mapsto B_v$ realizes an
intertwining between $\tau$ and $\Ind_{R_\t}^G(\o_\t).$
Conversely, if there is an embedding $T$ of $\tau$
into $\Ind_{R_\t}^G(\o_\t),$ then
setting $\l(v)=[T(v)](e),$ we get a
map $\l:V_{\tau}\rightarrow V_{\o}$
with the property specified above.
Thus, $\tau$ has an $\o_\t$\snug-Bessel model if
and only if a Bessel functional $\l$ exists.
\endroster
\endremark

In this paper we shall make use of the following basic uniqueness
principle.

\proclaim{Theorem/Conjecture 1.4}
Let $\tau\in\Cal E(G).$ 
Then for a fixed
$\o$ and $\t,$ we have
$\dim_\C\Hom_G(\tau,\Ind_{R_\t}^G(\o_\t))\leq 1.$
That is, a Bessel model is unique for irreducible representations.\qed
\endproclaim

Uniqueness for Whittaker models is well-known.  For rank
one Bessel models, Theorem 1.4 was proved, for both orthogonal
and unitary groups, by Novodvorsky\Lspace \Lcitemark Nov\Rcitemark \Rspace{}.
For Bessel models of arbitrary rank, Theorem 1.4 has been 
proved when $\G$ is
an orthogonal group by S. Rallis ({\rm\Lspace \Lcitemark Ral\Rcitemark \Rspace{}}).
Though the argument in the unitary case should be similar, it
has not yet been written down in full detail. 

{\it In the remainder of this paper we study those Bessel
models for which the uniqueness principle above is valid.
Thus we assume that Theorem/Conjecture 1.4 is true henceforth.}
Our results are therefore complete for split orthogonal groups
and for rank one Bessel models on quasi-split unitary
groups, while they are contingent upon the truth of
Theorem/Conjecture 1.4 for higher rank Bessel models
in the unitary case.

To conclude this section we introduce the notion of
a minimal Bessel model for an admissible representation
$\tau$ of $G$.  This will be a key notion in what follows.

\definition{Definition 1.5}  
Suppose that $\tau$ has an $\o_\chi$\snug-Bessel model which
is of rank $\ell_1$.
We say that this model is {\it minimal} if 
$\tau$ has no Bessel model of rank $\ell_1-1$
with respect to a representation $\o'_{\chi'}$ 
obtained as follows:  $\chi'$ is
a character of $U_{\ell+1}$ such that $\chi'=\chi$ on the
simple roots of $U_\ell$ (this implies that $M_{\chi'}\subset M_{\chi}$),
and $\o'$ is a component of $\o|_{M_{\chi'}}$.
\enddefinition

This condition is used in our proof of Proposition 2.4 below;
see the discussion following the proof of Lemma 2.10.

If $\tau$ has a Bessel model, we denote by $\Cal B(\tau)$ 
the smallest non-negative integer $\ell_1$ such that $\tau$
has a Bessel model of rank $\ell_1$.
For example, $\tau$ is generic if and only if 
$\Cal B(\tau)=0$.  
Then any Bessel model for $\tau$ of rank $\Cal B(\tau)$
is clearly a minimal Bessel model in the sense
of Definition 1.5.  In particular, any representation which
has a Bessel model has a minimal Bessel model.

\subheading{\Sec 2 Induction of Bessel models}
In this section we study the behavior of minimal Bessel models
under induction and prove an analogue of Rodier's Theorem\Lspace \Lcitemark Rodb\Rcitemark \Rspace{}
for such models.

Suppose that $\PP=\M\N$
is an arbitrary parabolic subgroup of $\G.$
Then 
$$\gather
\M\simeq GL_{n_1}\times\dots\times GL_{n_t}\times \G(m)\tag 2.1\\
\intertext{if $\G$ is orthogonal, and}
\M\simeq \Res^E_F(GL_{n_1})\times\dots\times \Res^E_F(GL_{n_t})\times \G(m)\tag 2.2
\endgather$$
if $\G$ is unitary,
where 
$$\G(m)=\cases SO_{2m+1}&\text{if }\G=SO_{2r+1};\\ SO_{2m}&\text{if }\G=SO_{2r};\\
U_{2m+1}&\text{if }\G=U_{2r+1};\\ U_{2m}&\text{if }\G=U_{2r},\endcases$$
and we take the convention that $SO_1=\{1\}.$
Here $r=n_1+\dots+n_t+m.$
Let $\pi\in\Cal E(M).$  Then
$$\pi=\s_1\otimes\dots\otimes\s_t\otimes\tau,\tag 2.3$$
where $\s_i\in\Cal E(GL_{n_i}(F))$ or $\Cal E(GL_{n_i}(E)),$ accordingly,
and $\tau\in\Cal E(G(m)).$  
Suppose that $\tau$ has a Bessel model.  We
let $\ell_1$ be the rank of a minimal Bessel
model for $\tau$, 
$\ell_0=2\ell_1+r_0-2r,$ and $\ell'=m-\ell_1.$
Let $\B'=\TT'\U'=\B\cap\G(m),$ and
$\U'_{\ell'}$ be the subgroup of $\U'$
consisting of matrices whose middle $\ell_0\times\ell_0$
block is the identity.
Choose a character $\t_1$ of $U'_{\ell'}$ and
$\o\in\Cal E(M_{\t_1})$ for which $\tau$ has
an $\o_{\t_1}$\snug-Bessel model which is minimal.
Suppose that each of the representations $\s_i$
is generic.  Let $\ell=r-\ell_1,$
and let $\t$ be a character of $U_\ell$
of the form $\t=\t_0\otimes\t_1,$
where $\t_0$ is a generic character on each $GL$ block
corresponding to a fixed non-trivial additive character $\psi$
of $F.$  (We call this the
$\psi$--{\it generic} character of the $GL$ component.)
Let $\tilde w_0$ be the longest element of $W(\G,\A_0)/W(\M,\A_0)$
and fix a representative $w_0$ for $\tilde w_0.$
Our first main result is the following.

\proclaim{Theorem 2.1}  Let $\eurm k=F$
if $\G$ is orthogonal and $E$ if $\G$ is unitary.
Let $\PP=\M\N$ be a parabolic subgroup of $\G,$
with $\M$ as in {\rm (2.1)} or {\rm (2.2)}.  Let $\pi$
be as in {\rm (2.3)} with each $\s_i$ generic.
Further suppose that $\tau$ has a Bessel model,
and that $\t_1$
is a character of $U_\ell\cap G(m)$ which
gives rise to an $\o_{\t_1}$--Bessel model for $\tau$
which is minimal.
Let $\t$ be a character of $U_\ell$ such that
$\t|_{U_\ell\cap GL_{n_i}(\eurm k)}$ is $\psi$--generic for each $i$ and such that
$\t|_{U_{\ell}\cap G(m)}=\t_1.$
Then $\Ind_P^G(\pi)$ has a unique $\o_\t^{w_0}$\snug-Bessel model.
Conversely, if any of the $\s_i$ are non-generic, or if $\tau$ 
has no Bessel model, 
then $\Ind_P^G(\pi)$ has no Bessel model.
\endproclaim

The remainder of this section will be devoted to the proof
of Theorem 2.1.  The first step is to 
reduce the theorem to the case of a  maximal proper parabolic
subgroup.  To do this, suppose Theorem 2.1
holds for maximal proper parabolic subgroups and let $\PP=\M\N$
be an arbitrary parabolic.  Then $\M$ is of the form (2.1) or (2.2).
Let $\PP_1=\M_1\N_1$ be the standard maximal proper parabolic with
$\M_1=GL_{r-m}\times \G(m)$ or $\M_1=\Res^E_F(GL_{r-m})\times\G(m)$
which contains $\M.$
Let $\rho=\Ind_{P\cap M_1}^{M_1}(\pi).$  Then $\rho=\rho_1\otimes\tau,$ 
where $\rho_1$ is the representation of $GL_{r-m}(\eurm k)$
parabolically induced from $\s_1\otimes\dots\otimes\s_t.$  Since each $\s_i$ is generic, Rodier's Theorem
implies that $\rho_1$ has a unique generic constituent.
Now for each irreducible constituent $\pi_1$ of $\rho_1,$
the representation $\pi_1\otimes\tau$ satisfies the hypothesis
of the Theorem.
Then, by assumption, 
$$\Ind_P^G(\pi)=\Ind_{P_1}^G(\Ind_{P\cap M_1}^{M_1}(\pi)\otimes\bold 1_{N_1})$$
will have a unique Bessel model of the desired type.
%Also, if $\Ind_P^G(\pi)$ had a model of rank
%less than $\Cal B(\tau)$, then the maximal parabolic
%case would imply that $\Ind_{P\cap M_1}^{M_1}(\pi)\otimes
%\bold 1_{N_1}$ did, and hence $\tau$ did too, a
%contradiction.  Thus the model so-obtained is minimal.

Now suppose that
$\PP=\M\N$ is a maximal proper parabolic
subgroup of $\G.$
Then for some  $n,$ $1\leq n\leq r,$ and
$m=r-n$ we have $\M\simeq GL_n\times \G(m)$ if
$\G$ is orthogonal, and $\M\simeq\text{Res}_{E/F}(GL_n)\times \G(m)$
if $\G$ is unitary.  
Let $\pi=\s\otimes\tau\in\Cal E(M),$ where $\s\in \Cal E(GL_n(\eurm k))$
and $\tau\in\Cal E(G(m)).$  
Suppose that
$\tau$ has an $\o_{\t_1}$\snug-Bessel model
of rank $\ell_1\geq0,$ and it is minimal.
Assume that $\t_1$ is of the form given in Lemma
1.2.
Let $\ell=r-\ell_1.$
Note that $\ell_1\leq m$ implies $\ell\geq n.$
Let $\ell'=\ell-n=m-\ell_1.$  Then $\t_1$ is a character
of $U'_{\ell'},$ where $\U'_{\ell'}=\U_{\ell'}\cap \G(m).$
Let $\t_0$ be the  generic character of 
the upper triangular unipotent subgroup $\U_0$ of $GL_n$
by given by a fixed additive character 
$\p.$    Now define
the character $\t$ on $U_\ell$ by
$\t=\t_0\otimes\t_1\otimes 1_{U'},$ where $U'$ is the
complement of
$U_0\times U'_{\ell'}$ in $U_\ell.$
Note that $M_\t=M_{\t_1}.$
We will examine the space of $\o_\t$\snug-Bessel functionals
for $\Ind_P^G(\s\otimes\tau).$

In order to carry out our computation, we have
to give a description of the $R_\t-P$ double cosets
in $G.$  The given  list  of such double cosets is
exhaustive, but overdefined.

Let $W_{\M}=W(\M,\TT_d).$  Then
$$W_{\M}\simeq\cases S_n\times (S_m\ltimes \Z_2^m)&\text{if }\G\neq SO_{2r}\\
S_n\times(S_m\ltimes\Z_2^{m-1})&\text{otherwise.}\endcases$$
Note that $|W/W_{\M}|=2^n\pmatrix r\\n\endpmatrix.$
Let $w_0$ denote the longest element of $W/W_{\M}.$
Then 
$$w_0=(1\,r_0)(2\ r_0-1)\dots (n\, r_0+1-n),$$
unless $\G=SO_{2r}$ and $n$ is odd, in which case
$$w_0=(1\,r_0)(2\ r_0-1)\dots (n\, r_0+1-n)(r\, r+1).$$
We now give a list of coset representatives for $W/W_{\M}.$
We will say that a permutation $s\in S_{r_0}$ ``appears''
in $w$ if $w=w's,$ for some $w'$ which is disjoint from $s.$
We will also use the convention that if $1\leq i\leq r,$ then
$i'=r_0+1-i.$

\proclaim{Lemma 2.2} Suppose that $w\in W.$  
\roster
\item"(a)"  If $\G\neq SO_{2r},$ then  there is
an element $w_1$ of $W$ so that $w\equiv w_1\, (\mod W_{\M})$
with $w_1$ a product of disjoint transpositions in $S_{r_0}.$
More precisely,
we may choose $w_1$ of the form $w_1'w_1'',$
with
$$w_1'=\prod_{i=1}^k (a_i\ a'_i),$$
for some $\{a_i\}\subset \{1,\dots,n\},$
and 
$$w_1''=\prod_{i=1}^j(b_i\, c_i)(c'_i\, b'_i),$$
with $\{b_i\}\subset\{1,\dots,n\},$ and 
$\{c_i\}\subset\{n+1,n+2,\dots,r_0-n\}.$
Furthermore we may assume that the transpositions appearing in $w_1'$
and $w_1''$ are all disjoint.
\item"(b)" If $\G=SO_{2r},$ then $w\equiv w_1w_2,$
where $w_1$ is of the form given in part (a),
and either $w_2=1,$ $w_2=(d_0\,d_0'),$ for some $n+1\leq d_0\leq r,$
or $w_2=(i_0\,j_0'\,i_0'\,j_0),$ for some
$1\leq i_0\leq n<j_0\leq r.$  In each case $w_1$ and $w_2$
are disjoint.
\endroster
\endproclaim

\demo{Proof}
We first write $w=cs,$ with $s\in S_r$ and $c\in \Z_2^r.$
Since $c$ acts on the cycles of $s$ independently,
we may assume that $s$ is a pair of
``companion'' cycles, $(a_1\,a_2\dots a_t)(a_1'\,a_2'\dots a_t').$
If $s=1,$ or the length of each of the
two companion cycles in $s$ is two, then the claim
is trivially true, so we assume that the length
of each of the cycles is greater than two.
Suppose that the claim holds whenever
the length of the two cycles in $s$ is less
than $t.$
Without loss of generality, we may assume that $a_1\leq n.$
If, for some $i,$ we have $a_i,a_{i+1}\leq n,$ then
$$w\equiv w(a_i\, a_{i+1})(a_i'\,a_{i+1}')=c( a_1\dots a_{i-1}\, a_i\, a_{i+2}\dots a_t)( a_1'\dots a_{i-1}'\, a_i'\, a_{i+2}'\dots a_t'),$$
and the claim holds by induction.  Similarly,
we may assume that if $ a_i>n,$ then $ a_{i+1}\leq n.$
This argument also shows that we may assume that $t$
is even.
Now we see that
$$\align
w&\equiv cs\cdot ( a_1\,a_{t-1}\, a_{t-3}\dots a_3)(a_1'\,a_{t-1}'\,a_{t-3}'\dots a_3')\\&=c( a_1\, a_t)( a_3\, a_2)\dots ( a_{t-1}\, a_{t-2})( a_1'\, a_t')( a_3'\, a_2')\dots ( a_{t-1}'\, a_{t-2}').
\endalign$$
Now  write $c=(b_1\,b_1')(b_2\,b_2')\dots (b_s\, b_s'),$
with $b_i\neq b_j,$ for $i\neq j.$

If, for a fixed even $i\geq 2,$
$\{a_i,a_{i+1}\}\subset \{b_j\}_{j=1}^s,$
then the product
$$(a_i\, a_i')(a_{i+1}\,a_{i+1}')(a_{i+1}\,a_{i})(a_{i+1}'\,a_{i}')=(a_i\,a_{i+1}')(a_{i+1}\,a_{i}')$$
appears in the reduced product for $w.$
The same is true if $\{a_1,a_t\}\subset\{b_j\}_{j=1}^s,$
i.e.,
$(a_1\,a_t')(a_t\,a_1')$ appears in $w.$
If $i\geq 2$ is even and $\{a_i,a_{i+1}\}\cap\{b_j\}_{j=1}^s=\emptyset,$ then $c$ commutes with $(a_i\,a_{i+1})(a_i'\,a_{i+1}'),$ and so this product of transpositions appears in $w.$
Similarly, if $\{a_1,\,a_t\}\cap\{b_j\}_{j=1}^s=\emptyset,$
then $(a_1\,a_t)(a_1'\,a_t')$ appears in $w.$

Suppose $i\geq 2$ is even and
that exactly one element of $\{a_i,a_{i+1}\}$ belongs
to $\{b_j\}_{j=1}^s.$  
Then, if $\G\neq SO_{2r},$ we can replace $w$
by $w(a_i\, a_i'),$ and we see that either $(a_{i+1}\,a_i)(a_i'\,a_{i+1}')$
or $(a_{i+1}\,a_i')(a_i\, a_{i+1}')$ appears in $w(a_i\, a_i'),$
depending on whether $a_{i+1}$ or $a_i$ is in $\{b_j\}_{j=1}^s.$
If $\G=SO_{2r}$ and $w$ either fixes some $d_0>n,$
or interchanges some $d_0$ and $d_0',$
then we can instead multiply $w$ by $(a_i\, a_i')(d_0\, d_0'),$
which shows that one of $(a_{i+1}\,a_i)(a_i'\,a_{i+1}')$
or $(a_{i+1}\,a_i')(a_i\, a_{i+1}')$ appears.
We see that the above considerations apply equally well to the pair
$\{a_1,a_t\}.$
By fixing the element $d_0$ before starting the above process, we can
guarantee that, when we have concluded, $w_2=1$ or $w_2=(d_0\, d_0').$

Finally suppose that no such $d_0$ exists.
Thus, $w(d)\neq d,d'$ for all $n+1\leq d\leq r.$  
So we may now assume that $d\in\{a_i\},$ for each $d,$
$n+1\leq d\leq r.$
Suppose that the number of $d$ for which $d=a_i\in \{b_j\}$ 
with $a_{i+1}\not\in \{b_j\}$ is even. (Here we are
including $\{a_1,\, a_t\}$  as one possible pair.)
Then we see that
$$w\equiv w\prod (d\, d'),$$
where the product is over precisely those $d=a_i$
for which $a_{i+1}\not\in \{b_j\},$
is of the form $w_1$ as claimed.
Finally if 
$\big|\{n+1\leq d\leq r|d\in \{b_j\}\}\big|$ is odd,
then we fix some such $d_0.$  Without loss of generality, assume that
$d_0=a_t.$ 
Multiplying on the right by the elements $(d\, d'),$
for the other such $d,$ we see that we have a factor
of $(a_t\, a_t')(a_1\, a_t)(a_1'\, a_t')$ remaining to be dealt with.
But this product is indeed $w_2=(a_1\, a_t'\, a_1'\, a_t),$ as claimed.
\qed
\enddemo

If $\G=SO_{2r},$ and 
$w_2$ is of this final form, then there is some flexibility
as to the indices appearing in $w_2.$
That is, we may choose, for $d_0,$ any
of the $a_i>n$ for which
$(a_i\, a_i')$ appears in $c,$
but $(a_{i+1}\, a_{i+1}')$
does not. We will need this below.

Recall that $\ell_0=r_0-2\ell.$  Let $s=\ell_0-2.$
Suppose $\x=(x_1,\dots,x_{s})\in F^{s}$
if $\G$ is orthogonal and $\x\in E^s$ if $\G$
is unitary.
Let 
$$n(\x)=\pmatrix I_{\ell}\\
&1&x_1&\dots&x_s&*\\
&&1&0&\dots&-\bar x_s\\
&&&\ddots&0&\vdots\\
&&&&1&-\bar x_1\\
&&&&&1\\
&&&&&&I_{\ell}\endpmatrix,$$
where $\bar x$ is the Galois conjugate of $x$
if $\G$ is unitary, and is $x$ if $\G$ is orthogonal.

The next result is a straightforward consequence
of Witt's Theorem and the Bruhat decomposition.
Here and for the rest of this section, we pass between
a Weyl group element and its coset representative without
changing the notation.

\proclaim{Proposition 2.3}
\roster
\item"(a)"  Let $g\in G.$   Then for some $w\in W$
and some $\x\in F^s,$ we have $R_\t gP=R_\t n(\x)wP.$
Clearly we can choose $w$ up to $W_{\M},$
i.e., we may assume $w=w_1$ is of
the form given in Lemma 2.2.
\item"(b)"  
Denote by $||\x||$ the standard length of
$\x\in F^s$ or $E^s$ accordingly.
If $||\x||=||\x_1||,$ then
$R_\t n(\x)wP=R_\t n(\x_1)wP,$ for all $w.$
\endroster
\endproclaim

If $H\subset G,$ we will use $\ind_H^G(\pi)$ to denote the representation
of $G$ compactly induced from $\pi$\Lspace \Lcitemark BeZa\Citecomma
Cas\Rcitemark \Rspace{}.
Recall that $\Ind_P^G(\pi)=\ind_P^G(\pi),$
by the Iwasawa decomposition.
If $V$ is a complex vector space, let $C^\infty(G,V)$
denote the space of locally constant $V$--valued functions on $G,$
and let  $C^\infty_c(G,V)$  
denote the subspace of elements of
$C^\infty(G,V)$ with compact support.  Let $\DD(G,V)=C^{\infty}_c(G,V)^*$
be the space of $V$--distributions on $G.$

Let $V_\s$ be the space of $\s,$  $V_\tau$ be the space of
$\tau,$ and $V_{\o}$ the space of $\o$
(and hence the space of $\o_\t$).
We let $V_{\pi}=V_\s\otimes V_\tau.$
Denote by
$V$ the vector space $\tilde V_{\o}\otimes V_{\pi},$
where $\tilde V_\o$ is the space of the smooth contragredient $\tilde \o$
of $\o.$

We wish to analyze the space
$\Hom_G(\Ind_P^G(\pi),\Ind_{R_\t}^G(\o_\t)).$
Dualizing, and using Theorem 2.4.2 of\Lspace \Lcitemark Cas\Rcitemark \Rspace{},
this is isomorphic to the space
$$\Hom_G(\ind_{R_\t}^G(\tilde\o_\t),\Ind_P^G(\tilde\s\otimes\tilde\tau)).$$
This space, in turn, is isomorphic to the space of intertwining forms on
$$\ind_{R_\t}^G(\tilde\o_\t) \bigotimes\Ind_P^G(\pi)$$
\Lcitemark Har\LIcitemark{}, Lemma 4\RIcitemark \Rcitemark \Rspace{}.
Now by Bruhat's thesis (see\Lspace \Lcitemark Rodb\LIcitemark{}, Theorem 4\RIcitemark \Rcitemark \Rspace{})
this is isomorphic to the space of $V$--distributions $T$ on $G$
satisfying
$$\e(r)*T*\e(p^{-1})=\d_{P}^{1/2}(p)T\circ[\tilde\o_\t(r)\otimes\pi(p)],\tag 2.4$$
for all $r\in R_\t$ and $p\in P.$

The analysis of this space of distributions will
make use of the following proposition.  Its proof requires
a combinatorial argument, and will be given
in several steps
later in this section.
\proclaim{Proposition 2.4} 
If there is a non-zero $V$\snug-distribution
$T$ satisfying {\rm (2.4)}
for all $r\in R_\t$ and $p\in P$  
which is supported on $R_\t n(\x) w P,$ then
$R_\t n(\x)wP=R_{\t} w_0 P$ and $\s$ is generic.
\endproclaim

\proclaim {Lemma 2.5}  Suppose that $T$ satisfies {\rm (2.4)}.
Then $T$ is completely determined by
its restriction to $R_\t w_0P.$
\endproclaim
\demo{Proof}  First note that 
a straightforward matrix computation shows that
$n(\x)w_0=w_0n(\x),$ for any $\x.$
Thus 
$$R_\t w_0P=\bigcup_{\x}R_{\t}n(\x)w_0 P=Pw_0P,$$
is open.  Therefore 
$C=G\setminus R_\t {w_0}P$ is closed.  Therefore, we have the
exact sequence\Lspace \Lcitemark BeZa\LIcitemark{}, \Sec 1.7\RIcitemark \Rcitemark \Rspace{}
$$0\longrightarrow C_c^\infty(R_\t{w_0}P)\longrightarrow C_c^\infty(G)\longrightarrow C_c^\infty(C)\longrightarrow 0.$$
Then, by tensoring with $V,$ the above exact sequence yields
the exact sequence
$$0\longrightarrow C_c^\infty(R_\t {w_0P},V)\longrightarrow C_c^\infty(G,V)\longrightarrow C_c^\infty(C,V)\longrightarrow 0.$$
Dualizing, we get the exact sequence
$$0\longrightarrow\DD(C,V)\longrightarrow\DD(G,V)\longrightarrow\DD(R_\t w_0 P,V)\longrightarrow 0.$$
Let $\DD_{R_\t,P}$ be  the subspace of distributions
satisfying (2.4).  Then Proposition 2.4
implies that if $T\in\DD(G,V)_{R_\t,P}$ and $f\in C_c^\infty(C,V),$
then $T(f)=0.$  Thus,
the above sequence tells us that $\DD(G,V)_{R_\t,P}\simeq\DD(R_\t w_0P,V)_{R_\t,P},$
which completes the proof of the Lemma.\qed
\enddemo

Let $R_\t^{w_0}=w_0^{-1}R_\t w_0,$ and
denote by $\o_\t^{w_0}$ the representation of
$R_\t^{w_0}$ defined by $\o_\t^{w_0}(r)=\o_\t(w_0rw_0^{-1}).$
Recall that $\bold P=\M\N$ is the
Levi decomposition of $\bold P.$

\proclaim{Lemma 2.6}  There exists an isomorphism between 
the vector space $\DD(R_\t w_0P,V)_{R_\t,P}$
and the vector space of distributions in
$\DD(U_\ell)\otimes\DD(P,V)$ of the form
$$\t(u)\, du\otimes\d_{P}^{-1/2}(m)\, dQ(m)\, dn,$$
where $Q\in\DD(M,V)$ satisfies
$$\e(r)*Q*\e(m^{-1})=Q\circ[\tilde\o^{w_0}_\t(r)\otimes\pi(m)],\tag 2.5$$
for all $r\in R_\t^{w_0}\cap M,$ $m\in M.$
\endproclaim

\demo{Proof}  Define a projection
$$\PPP:\cc(U_\ell)\otimes \cc(P,V)\to\cc(U_\ell w_0P,V)$$
by specifying that for all $f_1\in\cc(U_\ell)$
and $f_2\in\cc(P,V)$, one has
$$\PPP(f_1\otimes f_2)(uw_0p)=
\int\limits_{U_\ell\cap w_0Pw_0^{-1}}f_1(uu_1)\,f_2(w_0^{-1}u_1^{-1}w_0p)
\,du_1.$$
Then it follows from\Lspace \Lcitemark Sil\LIcitemark{}, Lemma 1.2.1\RIcitemark \Rcitemark \Rspace{} that $\PPP$ is
onto.  Let $T\in \DD(R_\t w_0P,V)_{R_\t ,P}$.
For $f_1, f_2$ as above, define $T'\in\DD(U_\ell)\otimes\DD(P,V)$ by
$$T'(f_1\otimes f_2)=T(\PPP(f_1\otimes f_2)).$$
Then one sees easily that (2.4) implies the equality
$$\e(u)*T'*\e(p^{-1})
=\tot(u)T \circ[\pi(p)]\tag 2.6$$
for all $u\in U_\ell$, $p\in P$ (where $\pi$ acts on the second
factor of $V$).  As in\Lspace \Lcitemark Sil\LIcitemark{}, Section 1.8\RIcitemark \Rcitemark \Rspace{}, this implies that
$T'$ is in fact a pure tensor of the form
$$\chi(u)\,du\otimes\delta_P(m)^{-1/2}\,dQ(m)\,dn.\tag2.7$$
where $Q\in \DD(M,V)$.  (Here we are using that
$\pi(mn)=\pi(m)$.)  It is a formal consequence of the
definitions that (2.6) implies that
$$Q*\e(m^{-1})=Q\circ[\pi(m)]$$
for all $m\in M$.  We claim that, more strongly, equation
(2.5) holds.  To see this, write
$$dQ(p)=\delta_{P}(m)^{-1/2}\,dQ(m)\,dn.$$
Let $f_1\in \cc(U_\ell)$, $f_2\in \cc(P,V)$ and $r\in R_\t\cap w_0Mw_0^{-1}.$
Then by (2.4) we have
$$\align
\int\limits_{U_\ell}f_1(u)\t(u)\,du\int\limits_P\tilde\o_\t(r)f_2(p)\,dQ(p)&=
\int\limits_{U_\ell\times P}f_1(u)\, \tilde\o_\t(r)f_2(p)\,dT'(u,p)\\
&= T\left(\tilde\o_\t(r)\PPP(f_1\otimes f_2)\right)\\
&= \int\limits_{U_\ell w_0P}\PPP(f_1\otimes f_2)(ruw_0p)\,dT(uw_0p).
\endalign$$
But $ruw_0p=(rur^{-1})w_0(w_0^{-1}rw_0p)$, so this
expression is equal to
$$\align
\int\limits_{U_\ell\times P}f_1(rur^{-1})f_2(w_0^{-1}rw_0p)&\,dT'(u,p)\\
&=\int\limits_{U_\ell}f_1(rur^{-1})\t(u)\,du\,\int\limits_Pf_2(w_0^{-1}rw_0p)\,dQ(p)\\
&=\int\limits_{U_\ell}f_1(u)\t(u)\,du\int\limits_P f_2(p)\,d(\e(w_0^{-1}rw_0)*Q)(p),
\endalign$$
where in this last equality the defining properties of $R_\t=M_\t U_\ell$
have been used to simplify the $U_\ell$ integral.
Since this holds for all $f_1\in\cc(U_\ell)$ one concludes that
$$\e(w_0^{-1}rw_0)*Q=Q\circ\tilde\o_\t(r)$$
for all $r\in R_\t\cap w_0Mw_0^{-1}$, as desired.

Conversely, given a distribution $Q$ satisfying equation (2.5),
one reverses the above steps to arrive at a 
distribution $T'\in\DD(U_\ell)\otimes\DD(P,V)$ satisfying (2.6).
Since the map $\PPP$ is onto, one may define a distribution
$T\in\DD(U_\ell w_0P,V)$ by the formula
$$T(\PPP(f_1\otimes f_2))=T'(f_1\otimes f_2)$$
provided one shows that if $\PPP(\sum_if_{1,i}\otimes f_{2,i})=0$,
then $T'(\sum_if_{1,i}\otimes f_{2,i})=0$.  This follows as
in\Lspace \Lcitemark HeR\LIcitemark{}, Theorem 15.24\RIcitemark \Rcitemark \Rspace{}.
Since $M_\t\subseteq w_0^{-1}Mw_0$ and $R_\t=U_\ell M_\t$,
it follows from (2.5) and (2.7) that the $T$ so-obtained
satisfies (2.4).

The maps $T\mapsto Q$, $Q\mapsto T$ described above are clearly inverses.
This completes the proof of the Lemma.\qed\enddemo

We now complete the proof of Theorem 2.1, modulo the proof of Proposition 2.4.  
Let $Q$
be as in the proof of Lemma 2.6.  Then by
Bruhat's thesis once again,  $Q$ corresponds to an element of
$$\Hom_{M}(\ind_{R_\t^{w_0}\cap M}^M(\tilde \o_\t^{w_0}),\tilde\s\otimes\tilde\tau),$$
which, by duality gives an element of
$\Hom_{M}(\pi,\Ind_{R_\t^{w_0}\cap M}^M(\o_\t^{w_0})).$
Since $M=GL_n(\eurm k)\times G(m),$ where $G_1$ is either $GL_n(F)$
or $GL_n(E),$ depending on whether $\G$ is orthogonal
or unitary, we see that
this last space is exactly the space of Whittaker models for $\s$
tensored with the space of $\o_\t^{w_0}$--Bessel models for $\tau.$
%The minimality of the $\omega_\chi^{w_0}$-Bessel
%model for $\Ind_P^G(\pi)$ follows from this argument
%and the proof
%of Proposition 2.4 below.
\qed

\demo{Proof of Proposition 2.4}
The remainder of the section will consist of
a proof of Proposition 2.4.
This is carried out in several steps.
We begin by showing that, 
on many double
cosets, the compatibility condition
$\pi(p)=\o_\t(wpw^{-1})$ can not
be satisfied for some $p\in P$ with $r=wpw^{-1}\in R_\t.$
By
\Lcitemark Sil\LIcitemark{}, Theorem 1.9.5\RIcitemark \Rcitemark \Rspace{}, this is sufficient to imply
the Proposition.

Let $\S_{\PP}^+$ denote the set of positive roots
in $\N.$  Let $\D$ denote the simple roots
of $\TT$ in $\G$ which give rise to our
choice of Borel subgroup. 
If $\a\in \P(\G,\TT),$ then we let $X_\a$
be the corresponding element of a Chevalley basis for the Lie algebra
of $\U$ or $\bar{\U}.$  
Let $\a_i$ denote the root $e_i-e_{i+1},$ and $\b=e_{\ell}+e_{\ell+1}.$
Let $X=\{\a_1,\a_2,\dots,\a_\ell,\b\}.$
Then $X$ is the set of roots where the character
$\t$ is non-trivial.
For $\a\in X$ we have
$\t(I+tX_\a)=\p_{\a}(t).$ 
Also, note that $X\cap\S_{\PP}^+=\{\a_n\}.$
We list the elements of $\S_{\PP}^+,$ for future reference.
If $\G=SO_{2r+1},$ then
$$\align
\S_{\PP}^+=&\left\{e_i\pm e_j\,|\,1\leq i\leq n<j\leq r\right\}\cup\\
           &\left\{e_i+e_j\, |\, 1\leq i<j\leq n\right\}\cup\left\{e_i|1\leq i\leq n\right\}.
\endalign$$
If $\G=U_{2r},$ then
$$\align
\S_{\PP}^+=&\left\{e_i\pm e_j\,|\,1\leq i\leq n<j\leq r\right\}\cup\\
           &\left\{e_i+e_j\, |\, 1\leq i<j\leq n\right\}\cup\left\{2e_i|1\leq i\leq n\right\}.
\endalign$$
If $\G=U_{2r+1},$ then
$$\align
\S_{\PP}^+=&\left\{e_i\pm e_j\,|\,1\leq i\leq n<j\leq r\right\}\cup\\
           &\left\{e_i+e_j\, |\, 1\leq i<j\leq n\right\}\cup\left\{e_i,2e_i|1\leq i\leq n\right\}.
\endalign$$
Finally, if $\G=SO_{2r},$
then 
$$\align
\S_{\PP}^+=&\left\{e_i\pm e_j\,|\,1\leq i\leq n<j\leq r\right\}\cup\\
           &\left\{e_i+e_j\, |\, 1\leq i<j\leq n\right\}.
\endalign$$
We list the various $I+X_\a,$ which generate the root subgroups
$\U_\a$ of $\U.$
Let $E_{ij}$ denote the elementary matrix
whose only non-zero entry is a $1$ in the $ij$\snug-th entry.
We recall the convention that $i'=r_0+1-i.$
Suppose $E=F(\g),$ where 
$\bar\g=-\g,$ and $a\mapsto\bar a$ is the
Galois automorphism of $E/F.$
If $\a=e_i-e_j,$ then $I+X_\a=I+E_{ij}-E_{j'i'}$ if $\G$ is
orthogonal, and $I+X_\a=I+\g E_{ij}-\g E_{j'i'}$
if $\G$ is unitary.
If $\a=e_i+e_j,$ then $I+X_\a=I+E_{ij'}-E_{ji'}$
if $\G$ is orthogonal, and $I+X_\a=I+\g E_{ij'}-\g E_{ji'}$
if $\G$ is unitary.
If  $\a=e_i,$
then $I+X_\a=I+E_{i,r+1}-E_{r+1,i'}$ if $\G=SO_{2r+1},$
and $I+X_\a=I+\g E_{i,r+1} -\g E_{r+1, i'}$
if $\G=U_{2r+1}.$
Finally, if $\G$ is unitary and $\a=2e_i,$
then $I+X_\a=I+\g E_{ii'}.$

Suppose that,
for some $\a\in\S_{\PP}^+,$ we have $\a'=w\a\in X.$
Choose some $t\in F^\times$ for which
$\p_{\a'}(t)\neq 1.$ 
Now set $p=I+tX_\a,$ which is in $P.$
Then $r=wpw^{-1}=I+tX_{\a'}\in U_\ell\subset R_\t.$
Note that $\pi(p)=1\neq\o_\t(r)=\p_\a(t).$
Thus, if $w$ has the above property,
$R_\t wP$ can support no
distribution of the desired type.

\proclaim{Lemma 2.7}  Let $\G=SO_{2r}.$
Suppose that,
as in Lemma 2.2, $w\in W$ is equivalent $\mod W_{\M}$ to 
$w_1w_2,$ with
$w_2=(i_0\, j_0'\, i_0'\, j_0),$ for some $1\leq i_0\leq n<j_0\leq r.$
Then $w\S_{\PP}^+\cap X\neq \emptyset.$
\endproclaim

\demo{Proof}From the proof of Lemma 2.2, we may assume that for
each $n+1\leq k\leq r,$ we have $w(k)=i_k$ or $i_k',$
for some $1\leq i_k\leq n.$  First suppose that
$(i_{n+1}\, n+1)(i_{n+1}'\, (n+1)')$ appears in $w.$
Consider first the case 
that for all $k,$ $n+1\leq k\leq \ell,$ we have
a permutation $(i_k\, k)(i_k'\,k')$ appearing in $w.$
Since $\ell+1=w(i_{\ell+1})$ or $\ell+1=w(i_{\ell+1}'),$
we have
$w(e_{i_{\ell}}+e_{i_{\ell+1}})=e_\ell\pm e_{\ell+1},$ which will be in
$X.$
So now we may suppose that either $j_0\leq \ell,$ or
$(i_k\, k')(i_k'\,k)$ appears in $w,$ for some $k$ with $n+1\leq k\leq \ell.$
Since $w$ changes an even number of signs,
we see that in the former case there must be some $k$
with
$n+1\leq k\leq r,$ so that $(i_k\, k')(i_k'\, k)$
appears in $w.$
Now we can multiply on the right by
$(j_0\, j_0')(k\, k'),$ to see that, in fact, we may assume
that $(i_{k_0}\, k_0')(i_{k_0}'\, k_0)$ is appearing, for some $k_0,$
with  $n+1\leq k_0\leq \ell.$
Choosing the minimal such $k_0,$ we know that
$k_0=w(i_{k_0}'),$ while $k_0-1=w(i_{k_0-1}).$  Thus,
$\a_{k_0-1}=w(e_{i_{k_0-1}}+e_{i_{k_0}})\in w\S_{\PP}^+\cap X,$ and the Lemma holds.

Thus, we may assume that either $(i_{n+1}\, (n+1)')(i_{n+1}'\, n+1)$
appears in $w,$ or that 
$(i_{n+1}\, (n+1)'\, i_{n+1}'\, n+1)$ does.
In the former case, we may multiply on the right
by $(n+1\, (n+1)')(j_0\, j_0'),$
to get an equivalent $w$ for which the latter is true, i.e,
we may assume that $i_0=i_{n+1}.$
First
suppose $w(n)=n.$  Then
$w(e_{n}+e_{i_0})=\a_n\in w\S_{\PP}^+\cap X$ and we are done.
Suppose instead  that $w(n)=n',$ i.e., that $(n\, n')$
appears in $w.$  Let $i$ be the smallest positive integer
so that $w(n-i)\neq n-i'.$ (By our assumption on the form of
$w,$ such an $i$ exists.) Then $n-i=i_k,$ for some $k\geq n+1,$
and  $n-i=w(k)$ or $w(k').$
Therefore,  $\a_{n-i}$ is equal to either $w(e_{n-i+1}-e_k)$
or to $w(e_{n-1+1}+e_k).$  In either case,
$\a_{n-i}\in w\S_{\PP}^+\cap X.$
Finally, we may suppose that either $n=i_0$
or that one of $(n\, k)(n'\, k')$ or $(n\, k')(n'\, k)$
appears in $w,$ for some $k,$ $n+1\leq k\leq r.$
If $(n\, k)(n'\, k')$ appears in $w,$
then we may multiply on the right by
$$((n+1)\, k)((n+1)'\, k')(i_0\, n)(i_0'\, n'),$$
to replace $w$ by an equivalent element with $i_0=n.$
Similarly, if
$(n\, k')(n'\, k)$ appears in $w,$
then we may multiply on the right by
$$(n+1\, k)((n+1)'\, k')(i_0\, n)(i_0'\, n')(k\, k')(n+1\, (n+1)'),$$
to see that we may assume that $i_0=n.$
We are thus reduced to the case 
where  $(n\, (n+1)'\, n'\, n+1)$ appears in $w.$
In this case,  $w(e_n+e_{n+1})=\a_n\in w\S_{\PP}^+\cap X.$
Thus, in all cases, the Lemma holds.\qed\enddemo

\remark{Remark}
For future use we make note of the following
fact.  If $w$ is as in Lemma 2.7, and if
$w^{-1}\a_{\ell}\in\S_{\PP}^+,$
then the proof of Lemma 2.7
shows that either $w^{-1}\a_{\ell}=e_i+e_j,$
for some $i,j\leq n,$ or
that $w\S_{\PP}^+\cap X\neq\{\a_\ell\}.$
\endremark

We now describe those
$w$ which have the property that $w\S_{\PP}^+\cap X=\emptyset.$
By Lemmas 2.2 and 2.7, we may assume that $w$
is a product of disjoint transpositions.

\proclaim{Lemma 2.8}
Suppose that $w\in W$ is a representative
for a class in $W/W_{\M},$ and $w$ is in the form specified by Lemma 2.2.
Further suppose that, for all $\a\in\S^+_{\PP},$ we have $w\a\not\in X.$
Then the following hold:
\roster
\item"(a)" For all $k$ with  $n+1\leq k\leq \ell,$ we have
$w(k)>n;$
\item"(b)"
For all $i$ with  $1\leq i\leq n,$
we have $w(i)\neq i.$
\endroster
\endproclaim
\demo{Proof}
(a)  First suppose that $w(\ell)\leq n.$
If $w(\ell +1)\leq n,$ then $w(\b)\in\S_{\PP}^+,$
contradicting our choice of $w.$
If $w(\ell+1)=\ell+1,$ then again $w(\b)\in\S_{\PP}^+.$
Finally, if $w(\ell+1)\geq n',$ then $w(\a_\ell)\in\S_{\PP}^+.$
So we must have $w(\ell)>n.$

Now suppose that for some $k,$
$n+1\leq k\leq \ell -1,$
we have
$w(k)\leq n.$  If $w(k+1)=k+1,$ or $w(k+1)\geq n',$
then $w(\a_k)\in\S_{\PP}^+,$ which is a contradiction.
Therefore, $w(k+1)\leq n.$  However, this implies, by induction,
that $w(\ell)\leq n,$ which we have already seen is impossible.
Therefore, $w(k)>n.$

(b) Suppose that $w(i)=i$ for some $i,$
$1\leq i\leq n-1.$
If $w(i+1)\neq i+1,$ then $w(i+1)>n,$ and so
$w(\a_i)\in\S_{\PP}^+.$  Since this contradicts
our choice of $w,$ we have $w(i+1)=i+1.$
We may thus suppose that $w$ fixes $n.$  
Now by part (a), we have $w(n+1)\geq n+1,$ and therefore,
$w\a_n\in\S^+_{\PP}.$  This again is a contradiction, so $w$ cannot fix $n.$ 
Therefore, $w$ fixes none of the integers $1,2,\dots,n.$
\qed\enddemo

\proclaim{Lemma 2.9}
Suppose that $w$ is as in Lemma 2.8 and assume that
$w(n)\neq n'.$  Then for $n+1\leq k\leq\ell,$
we have $w(k)=k.$
\endproclaim

\demo{Proof}
By Lemma 2.8(a) it is enough to show that it is impossible that
$w(k)\geq n'$ for any such $k.$  Suppose to the contrary that 
there is some $k,$ with $n+1\leq k\leq \ell,$ for which
$w(k)\geq n'.$
Then there some $i\leq n,$ for which $k=w(i').$
If $w(k-1)=k-1,$ then $\a_{k-1}=w(e_{i}+e_{k-1})\in w\S_{\PP}^+\cap X.$
Since this contradicts our choice of $w,$ we must have
$w(k-1)\geq n'.$ 
Therefore, by (downwards) induction, $w(n+1)\geq n'.$ 
Set $w(n+1)=i'.$
Since $w(n)\neq n,$ and, by assumption, $w(n)\neq n',$ 
either
$w(n)=k$ or $w(n)=k'$ for some $k,$
with $n+1\leq k\leq n'-1.$
Therefore, $\a_n=w(e_{i}+e_{k})$ or $\a_n=w(e_i-e_k).$
Either one of these possibilities contradicts our assumption on $w.$
Thus $w(n+1)<n',$
which then implies the result of the Lemma.\qed
\enddemo

\proclaim{Lemma 2.10} Suppose that $w$ is as in Lemma 2.8.
Suppose that there is
some $i,$
$2\leq i\leq n,$ for which $w(i)=i'.$
Then $w(i-1)=(i-1)'.$\endproclaim

\demo{Proof}  Suppose $w(i-1)\neq (i-1)'.$
By Lemma 2.8(b),
we can choose  $k,$ 
with $n+1\leq k\leq r$ so that $w(i-1)=k$ or $w(i-1)=k'.$
Now $\a_{i-1}=w(e_i- e_{k})$  or $\a_{i-1}=w(e_i+e_k).$
Since this contradicts our choice of $w$ we conclude $w(i-1)=(i-1)'.$
\qed
\enddemo

Thus, if
$w$ is chosen as in Lemma 2.2 with $w\S_{\PP}^+\cap X= \emptyset$
and $w(n)=n',$ then $w=w_0.$  
If $w\S_{\PP}^+\cap X=\emptyset,$
and $w(n)\neq n',$
then by Lemma 2.9,
$w(\ell)=\ell.$  If $w(\ell+1)\neq \ell +1,$
then either $\a_{\ell}$ or $\b$ would be
of the form $w\a$ for some $\a\in\S_{\PP}^+.$
Consequently, $w(\ell+1)=\ell+1,$
and therefore $w(\a)=\a$ for $\a\in\{\a_{n+1}\dots,\a_{\ell},\b\}.$
Thus, for some
$i_0<n,$ and some $\{a_i\}\subset \{\ell+2,\dots,(\ell+2)'\},$
$$w=(1\, r_0)(2\, r_0-1)\dots (i_0\, i_0')(i_0+1\, a_{i_0+1})((i_0+1)'\,a_{i_0+1}')\dots (n\, a_n)(n'\, a_n').\tag 2.8$$
Let $a=a_n$ if $a_n\leq r,$ and  $a=a_n'$ otherwise.
Then $w\a_n=\pm e_a-e_{n+1}.$
Let $X_1=\{\a_{n+1}\dots,\a_{\ell},\b,w\a_n\}.$
Note that 
$$X_1=w(\{\a_n,\a_{n+1}\dots,\a_{\ell},\b\}),$$
and is thus a linearly independent subset of
the root system $\P(\G(m),\TT'),$
where we recall that $\TT'=\TT\cap\G(m).$
We extend $X_1\setminus\{\b\}$ to a set of simple roots for
$\G(m).$ Set $\B'=\TT_1\U''$ to be the corresponding
Borel subgroup of $\G(m),$ and suppose that
$U_{\ell'+1}''$ is the subgroup of $\U''$
which is conjugate to $\U_{\ell'+1}'$ and generated
by the elements of $X_1.$  (Recall that $U'_{\ell'}$
is the subgroup supporting the character $\t_1$
which gives rise to the model for $\tau.$)
Now let $\t'$ be the character of $U''_{\ell'+1}$
so that $\t'(I+tX_{\a})=\p_\a(t),$ for $\a\in X_1\setminus\{\a_\ell\},$
and $\t'(I+tX_{\a_\ell})=\p_{\a_\ell}(\d t).$
Let $M_{\t'}$ be the corresponding normalizer in $M_{\ell'+1}.$
Note that $M_{\t'}\subset M_\t.$
Suppose that $m'\in M_{\t'}.$  
If the distribution $T$
satisfies (2.4), then $\e(m')*T=T\circ\o(m').$
So for some component $\o'$ of $\o|_{M_{\t'}},$ we have
$$\e(r)*T*\e(h)= T\circ[\tilde\o'_{\t'}(r)\otimes\tau
(h)],$$
for all $h\in \G(m),$ and $r\in R_{\t'}=M_{\t'}U''_{\ell'+1}.$
If $T$
is non-zero,  this now implies
that $\tau$ has a Bessel model with respect to
$U'',\, \t'\, $ and $\o'.$
However, since
$U''_{\ell+1}$ is isomorphic to $U_{\ell+1},$
this  is a rank $\ell_1-1$ Bessel model for $\tau.$
This contradicts the minimality of the
$\o_{\chi_1}$-Bessel model for $\tau$.
Hence, no such $T$ exists.

Note that this argument shows that 
$\Ind_P^G(\pi)$
cannot have any Bessel model of rank less than
$\Cal B(\tau)$ supported on $R_\chi w P$.

Finally, suppose that $w=w_0.$  Let $u\in \bar U\cap GL_n(F).$
Set $r=w_0^{-1}uw_0.$  Then $r\in U_{\ell},$ and $\t(r)=\t^{w_0}_0(u).$
Since $$\e(r)*T=T\circ[\t_0^{w_0}(u)]=T*\e(u)=T\circ[\s(u)],$$
we see that $\s$ must be generic if $T$ is non-zero
\Lcitemark Rodb\Rcitemark \Rspace{}.
This completes the proof of Proposition 2.4 for the cosets $R_\t w P,$
with $w\in W/W_{\M}.$

We now examine the double cosets represented by $n(\x)w,$ where
$\x=(x_1,\dots,x_s)$ is a vector.
Recall that
$$n(\x)=\pmatrix I_{\ell}\\
&1&x_1&\dots&x_s&*\\
&0&1&0&\dots&-\bar x_s\\
&&&\ddots&0\\
&0&0&\dots&1&-\bar x_1\\
&0&0&0&\dots&1\\
&&&&&&I_{\ell}
\endpmatrix,$$
where $\bar x$ is the Galois conjugate of $x$ if $\G$
is unitary, and $\bar x=x$ if $\G$ is orthogonal.

We assume that $w$ is of the form
given in Lemma 2.2.
First note that if $\a\in X\setminus\{\a_{\ell}\},$
then $n(\x)(I+tX_\a)n(\x)^{-1}=I+tX_\a.$
Suppose that $w\S_{\PP}^+\cap(X\setminus\{\a_{\ell}\})\neq\emptyset.$
Choose $\a'\in\S_{\PP}^+$ with $w\a'=\a\in X\setminus\{\a_{\ell}\},$
and $t$ for which $\p_\a(t)\neq 1.$
Setting $p=I+tX_{\a'},$ we have
$n(\x)wpw^{-1}n(\x)^{-1}=I+tX_\a\in R_\t.$  Furthermore
$\o_\t(x)=\p_\a(t)\neq 1,$ while $\pi(p)=1.$  Thus, $R_\t n(\x)wP$
supports no distributions satisfying (2.4).

Now suppose that $w\S_{\PP}^+\cap X=\{\a_{\ell}\}.$
First suppose that $w^{-1}\a_{\ell}=e_i+e_j,$
with $i,j\leq n.$ Without loss of generality, assume that
$w(i)=\ell,$ and $w(j)=(\ell+1)'.$
Suppose that $\ell+2\leq k\leq r.$
If $w(k)=i_k\leq n$ then $w^{-1}(e_\ell+e_k)=e_i+e_{i_k}\in\S_{\PP}^+.$
If instead $w(k)=i_k'$ for some $i_k\leq n,$ then $w^{-1}(e_{\ell}-e_k)=e_i+e_{i_k}.$
Finally, if $w(k)=k,$ then $w^{-1}(e_\ell\pm e_k)=e_i\pm e_k\in\S_{\PP}^+.$
Choose $s_0\leq s$ for which $x_{s_0}\neq0.$
Let $y=x_{s_0}.$
Choose $k_0$ with the property that
either $w^{-1}(e_\ell+e_{k_0})$ or $w^{-1}(e_\ell-e_{k_0})$
is an element of $\S_{\PP}^+.$  Denote the root
$e_{\ell}\pm e_{k_0}$ as $\a_0,$ with $\pm$
chosen so that $w^{-1}\a_0\in\S_{\PP}^+.$
We may also assume that
$X_{\a_0}$ has $-1$ as its $(r+\ell_1,\ell+s_0)$ entry
(see Lemma 2.3).
Now note that 
$$n(\x)(I+tX_{\a_0})n(\x)^{-1}=(I+tX_{\a_0})(I+ytX_{\b}).$$
Thus, if $\p_{\a_0}(yt)\neq 0,$
and $p=I+tX_{w^{-1}\a_0}\in N,$ then $\pi(p)=1,$
while $$\o_{\t}(n(\x)wpw^{-1}n(\x)^{-1})=\p_{\a_0}(yt)\neq 1.$$
Consequently, $R_\t n(\x)wP$ cannot support 
a $V$\snug-distribution of the desired form.

We are left with the cases
$w\S_{\PP}^+\cap X=\{\a_\ell\},$
but $w^{-1}\a_\ell\neq e_i+e_j$
for all $i,j\leq n,$ or $w\S_{\PP}^+\cap X=\emptyset.$
For the second of these two cases, the form of $w$
is given by (2.8). 
In order to complete the proof
we will determine the form of $w$
in the first case.  To do so we need a few lemmas.

\proclaim{Lemma 2.11} Suppose that $w\S_{\PP}^+\cap X=\{\a_{\ell}\},$
but $w^{-1}\a_{\ell}\neq e_i+e_j,$ for all $i,j\leq n.$
Then $w(\ell)=\ell.$
\endproclaim
\demo{Proof}
If $w^{-1}(\ell)=j'$ for some $j\leq n,$
then $w^{-1}\a_{\ell}\not\in\S_{\PP}^+,$ which is a contradiction.
Suppose $w^{-1}(\ell)=j\leq n.$  If $w(\ell+1)=\ell+1,$
then
$w^{-1}(\b)=e_j+e_{\ell+1}\in\S_{\PP}^+,$ contradicting
our choice of $w.$  If $w(\ell+1)=i\leq n,$
then $w^{-1}\a_{\ell}\not\in\S_{\PP}^+,$ which
also contradicts our choice of $w.$
Finally, if $w^{-1}(\ell+1)=i'$ for
some $i\leq n,$ then
$w^{-1}\a_{\ell}=e_j+e_{i},$ which is again a contradiction.
Thus, $w(\ell)=\ell.$\qed
\enddemo

\proclaim{Lemma 2.12}  If $w$ is as in Lemma 2.11, then
for $n+1\leq k\leq \ell-1,$ we have $w^{-1}(k)>n.$
\endproclaim
\demo{Proof}
Suppose that $w^{-1}(k)=j\leq n.$
If $w(k+1)=k+1,$ then $w^{-1}\a_k=e_j-e_{k+1}\in\S_{\PP}^+.$
If $w(k+1)=i'$ for some $i\leq n,$
then
$w^{-1}\a_k=e_j+e_i.$  Either case contradicts our
hypotheses.
Therefore $w^{-1}(k+1)\leq n.$  Now by induction,
$w^{-1}(\ell-1)\leq n.$  On the other
hand, by Lemma 2.11,  $w(\ell)=\ell.$  Therefore
$w^{-1}\a_{\ell-1}\in\S_{\PP}^+,$ contradicting our choice of $w.$
Consequently, $w^{-1}(k)>n.$\qed
\enddemo

\proclaim{Lemma 2.13} Suppose that $w$ is as in Lemma 2.11.
Then $w(n)\neq n.$
\endproclaim

\demo{Proof}  Suppose that $w(n)=n.$  If $w(n+1)=n+1,$
then $w$ fixes $\a_n,$ which is in
the intersection
of $X$ and $\S_{\PP}^+.$  If $w^{-1}(n+1)=j'$ for
some $j\leq n,$ then
$w^{-1}\a_n=e_n+e_j.$  Both of these possibilities contradict our
choice of $w.$  By Lemma 2.12, $w^{-1}(n+1)>n,$
and so these are the only two choices for $w(n+1).$
Since each leads to a contradiction,  $w(n)\neq n.$\qed
\enddemo

\proclaim{Lemma 2.14} Suppose that $w$
is as in Lemma 2.11.
\roster
\item"(a)" For all $i\leq n$ we have $w(i)\neq i.$
\item"(b)" If $w(i_0)=i_0',$ for some $i_0\leq n$ 
 then $w(i)=i'$
for all $i\leq i_0.$
\endroster
\endproclaim

\demo{Proof}
(a) Suppose that $w(i)=i$  for some $i\leq n.$     
Choose the maximal such $i.$
By Lemma 2.13, $i<n.$
Suppose that $w(i+1)=(i+1)'.$  Then $w^{-1}\a_i=e_i+e_{i+1}\in\S_{\PP}^+.$
Thus in this case we have a contradiction.
If $w(i+1)=k$ or $w(i+1)=k'$ for some $n+1\leq k\leq r,$
then $w^{-1}\a_i=e_i\pm e_k\in\S_{\PP}^+.$  
This is also a contradiction, and hence no $i\leq n$
can be fixed by $w.$

(b)  Suppose that $w(i)=i',$ for some $i\leq n.$
If $w(i-1)=k$ or $k'$ for some $n+1\leq k\leq r,$
then $w^{-1}\a_{i-1}=e_i\pm e_k\in\S_{\PP}^+.$
But by part (a), $w(i-1)\neq i-1,$
so the only remaining possibility is
$w(i-1)=(i-1)'.$  This gives the
claim by induction.
\qed
\enddemo

\proclaim{Corollary 2.15}
If $w$ is as in Lemma 2.11, then $w(n)=k_0$
or $w(n)=k_0'$ for some $n+1\leq k_0\leq r.$
\qed
\endproclaim

\proclaim{Lemma 2.16}  Suppose that $w$ is as in Lemma 2.11.
Then $w(k)=k$
 for all $k$
with $n+1\leq k\leq \ell-1.$ 
\endproclaim

\demo{Proof}
Suppose that $w^{-1}(n+1)=j'$ for some $j\leq n.$
Then, by Corollary 2.15, 
$w^{-1}\a_n=e_j\pm e_{k_0}\in\S_{\PP}^+,$ contradicting our choice
of $w.$  Thus, by Lemma 2.12, $w(n+1)=n+1.$

Now suppose  $w^{-1}(k)=j_k'$ for some $k$
with $n+2\leq k\leq \ell-1,$
and some $j_k\leq n.$
If $w(k-1)=k-1,$ then $w^{-1}\a_{k-1}\in\S_{\PP}^+,$
which is a contradiction.  Therefore,
by Lemma 2.12, $w^{-1}(k-1)=j_{k-1}',$ for some $j_{k-1}\leq n.$
By induction, this gives $w(n+1)\neq n+1,$
while we have just shown that $w(n+1)=n+1.$  Therefore,
$w(k)=k.$\qed
\enddemo

\proclaim{Lemma 2.17}  
Suppose that  $w$ is as in Lemma 2.11.
\roster
\item"(a)" Suppose that $\G\neq SO_{2r}.$
Then, for some $n_1,$ with
$0\leq n_1<n,$ and some 
$$\{k_j|1\leq j\leq n-n_1\}\subset\{\ell+1,\ell+2,\dots,(\ell+1)'\},$$
we have
$$\align
w=&(1\, r_0)(2\, r_0-1)\dots (n_1\, n'_1)
(n_1+1\, k_1)(k_1'\, (n_1+1)')\dots\\
&(n\,k_{n_2})(k'_{n_2}\, n').
\endalign$$
Here $n=n_1+n_2.$
Furthermore, $ k_j=(\ell+1)'$ for some $j.$
\item"(b)" If $\G=SO_{2r},$ and we write $w=w_1w_2$
as in Lemma 2.2, then $w_2=1$ or $w_2=(d\, d'),$ for some $\ell+2\leq d\leq r.$
Furthermore $w_1$ is of the form
$$\align
w_1=&(1\, r_0)(2\, r_0-1)\dots (n_1\, n'_1)
(n_1+1\, k_1)(k_1'\, (n_1+1)')\dots\\
&(n\,k_{n_2})(k'_{n_2}\, n'),
\endalign$$
with $n=n_1+n_2,$ and the integers
$k_j$ are as in part (a).
Moreover, $k_j= (\ell+1)'$ for some $j.$ 
\endroster
\endproclaim

\demo{Proof}
First note that if $\G=SO_{2r},$ and $w=w_1w_2,$
then Lemma 2.16 and the remark
following Lemma 2.7 imply
that $w_2$ is not of the form
$(i\, j'\, i'\, j),$ for some $1\leq i\leq n<j\leq r.$
Moreover, since $w^{-1}\a_{\ell}\in\S_{\PP}^+,$
Lemmas 2.16 and 2.11 imply that if $w_2=(d\, d'),$
then $\ell+2\leq d\leq r.$  If $\G\neq SO_{2r},$
let $w_2=1.$

By Lemma 2.14(a), $w(i)\neq i$ for all $i\leq n.$
By Lemma 2.11,
Corollary 2.15, and Lemma 2.16, $w(n)=k$ or $k',$
for some $\ell+1\leq k\leq r.$
Let $n_1$ be the largest nonnegative integer for
which $n_1<n$ and $w(n_1)=n_1'.$
If $n_1>0,$ then by Lemma 2.14(b) 
$w=(1\, r_0)(2\, r_0-1)\dots (n_1\, n_1')w_2w',$
where $w'(i)=i$ for all $i\leq n_1,$ and $w_2$ and $w'$ are disjoint.
Now  $w'(i)\neq i$
and $w'(i)\neq i'$ for $n_1+1\leq i\leq n,$  
and therefore $n+1\leq w'(i)\leq n'-1.$
However, by Lemma 2.16, $\ell+1\leq w'(i)\leq (\ell+1)'.$  Thus,
$$w'=(n_1+1\,k_1)(k_1'\, n'_1-1)\dots (n\, k_{n_2})(k_{n_2}'\,n'),$$
as claimed.
Finally,  Lemma 2.11
implies  $w(\ell+1)\neq\ell+1,$ and so
we must have
$(w')^{-1}(\ell+1)=w^{-1}(\ell+1)=j',$ for some $j\leq n.$
\qed
\enddemo

We now finish the proof of Proposition 2.4.  If $w=w_0,$ then $n(\x)w_0=w_0n(\x),$
and since $n(\x)\in P,$ we have $R_\t w_0P=R_\t n(\x)w_0P.$
If $w\S_{\PP}^+\cap X=\{\a_{\ell}\},$ or $w\S_{\PP}^+\cap X=\emptyset,$
then Lemma 2.16 and equation (2.8) show that $w(e_n+e_\ell)=e_\ell\pm e_k,$
for some $k,$ with
$\ell+2\leq k\leq r.$ Let $\a=w(e_n+e_\ell),$
and denote $I+tX_{e_n+e_\ell}$ by $p.$
As before, choose
$\x_0$ so that $R_\t n(\x_0)wP=R_\t n(\x)wP,$
and such that
$\x_0$ has a non-zero entry $y$ with
$\o_\t(n(\x)wpw^{-1}n(\x)^{-1})=\p_\a(yt).$
Note that $\pi(p)=1.$  
Choosing $t$ for which $\p_\a(yt)\neq 1,$
we see that $R_\t n(\x)wP$ cannot support a $V$\snug-distribution
of the desired form. \qed
\enddemo

From the argument above, it is apparent that
$\Ind_P^G(\pi)$ cannot have a Bessel model
of rank less than $\Cal B(\tau)$.  Hence we
obtain the following Corollary.

\proclaim{Corollary 2.18}
Let the notation be as in Theorem 2.1.  Suppose
that the $\o_\chi$\snug-Bessel model for $\tau$ is
of rank $\Cal B(\tau)$.  Then the 
$\omega_\chi^{w_0}$\snug-Bessel model for $\Ind_P^G(\pi)$
is also minimal, and of rank $\Cal B(\tau)$.
\endproclaim

The proof of Theorem 2.1 also gives the following result.

\proclaim{Corollary 2.19}  For any $\s,\, \tau,\,\ell,\, \t,\,$
and $\o,$ the support of the twisted Jacquet functor
$\pi_{U_\ell,\t}$ is a finite number of double cosets.
\qed
\endproclaim

\subheading{\Sec 3  Holomorphicity and local coefficients}

In this section we prove the holomorphicity of the Bessel functional and the existence
of a local coefficient.  To do so, we first
adapt the argument used by
Banks\Lspace \Lcitemark Ban\Rcitemark \Rspace{} to prove the holomorphicity of Whittaker functions
for metaplectic covers of $GL_n.$  
Banks's result is an extension of Bernstein's Theorem,
which establishes the meromorphicity under uniqueness
and regularity hypotheses.
We show that the desired regularity holds in the
case of Bessel functionals.
We then use an argument
similar to Harish-Chandra's and to 
Shahidi's in
the generic case to establish the existence of
the local coefficient under certain conditions 
(Theorem 3.8).  Corollary 3.9 shows that
the local coefficient factors in a manner analogous to the generic
case.

Let $\G$ be as in Section 1.
We use the conventions found in\Lspace \Lcitemark Cas\LIcitemark{}, \Sec 1\RIcitemark \Citecomma
Shaa\Rcitemark \Rspace{}
for subsets of simple roots, Weyl groups, and arbitrary parabolic
subgroups.
Suppose that $\Delta$ is the collection of simple roots
corresponding to our choice  of Borel subgroup.   Let $\theta\subset \Delta$
be a collection of simple roots and set $\PP=\PP_{\theta}$.
Then $\PP$ has Levi decomposition $\PP=\M_{\T}\N_{\T},$  with 
$$\M=\M_\T\simeq GL_{n_1}\times\dots\times GL_{n_k}\times\G(m),$$
for some $n_i,\, m$ such that $r=n_1+\dots+n_k+m.$
We  abbreviate this by writing $\M\simeq \G_1\times\G(m).$
We also write $\N=\N_\T.$

Let $\A=\A_\T$ be the split component of $\M.$
Denote by
$\frak a_{\C}^*=\left(\frak a_{\T}\right)^*_\C$ the complexified dual of the real Lie algebra of
$\A,$ $q_F$  the residual characteristic
of $F,$ and denote by $H_P$ the Harish-Chandra homomorphism
\Lcitemark Har\Citecomma
Shaa\Rcitemark \Rspace{}.
Suppose that $\s\in\Cal E(G_1)$ and $\tau\in\Cal E(G(m))$,
and let $\pi=\sigma\otimes\tau$. 
For $\nu\in\frak a_{\C}^*,$ 
let $I(\nu,\pi,\theta)$ denote the induced representation
$$\Ind_{P}^G(\pi\otimes q_{F}^{<\nu, H_P()>})$$
and let $V(\nu,\pi,\theta)$ denote the space of associated
functions.  
We also use $\Pi_\nu$
to denote the representation $I(\nu,\pi,\theta).$

Assume that $\s$ is generic and that $\tau$ has an $\o_{\t'}$--Bessel model 
which is minimal and 
of rank $\ell_0.$
Let $\t$ be the character of $U_\ell$ whose restriction
to $U_\ell\cap G(m)$ is $\t'$ and whose restriction to
$G_1(F)\cap U_\ell$ is a $\psi$--generic character
$\t_1.$ We will construct a non-zero functional
$\L_\t(\nu,\pi,\theta)$ on 
$X_\nu=I(\nu,\pi,\theta)\otimes\tilde V_{\o}$
so that, for a certain character $\d$ of $M_\t,$
$$\L_\t(\nu,\pi,\theta)(\Pi_\nu(mu)(f_\nu\otimes\tilde v))=
\d(m)\t(u)^{-1}\L_\t(\nu,\pi,\theta)(f_\nu\otimes\tilde \o(m^{-1})\tilde v),$$
for all choices of $f_\nu\otimes\tilde v\in X_\nu$ and $mu\in R_{\t}.$
Then we will show in Theorem 3.6
that the function $\nu\mapsto\L(\nu,\pi,\T)(x_\nu)$ is holomorphic, for a holomorphic section $\nu\mapsto x_\nu.$

Let  $K=G(\Cal O_F),$
where $\Cal O_F$ is the ring of integers in $F.$
Then $K$ is a good maximal compact subgroup of $G$\Lspace \Lcitemark Cas\Rcitemark \Rspace{}.
Let $K_m$ be the corresponding $m$\snug-th
principal congruence
subgroup.
Then each $K_m$ is normal in $K.$
Let $\Gamma_m$ be a complete set of coset representatives
for $P\cap K\backslash K/K_m.$
Note that $\Gamma_m$ is of finite cardinality.
Let 
$$Y=\left\{f\in C^\infty(K,V_{\pi})\,\big|\, f(pk)=\left(\pi\right)(p)f(k),\,\forall p\in P\cap K,\, k\in K\right\}.$$
Then $F\mapsto F|_K$ is a $K$\snug-isomorphism from $V(\nu,\pi,\theta)$
to $Y,$
by the Iwasawa decomposition of $G.$  
We will define a certain functional 
on $Y,$ and use this realization to
define an associated functional  on $X_\nu.$
Let 
$$Y_m=\left\{f\in Y|f(kk_1)=f(k),\forall k\in K,k_1\in K_m\right\}.$$
Thus, $Y_m$ is the set of $K_m$--fixed vectors of $Y$ under the
action of $K.$
Furthermore, the Iwasawa
decomposition allows us to realize
$I(\nu,\pi,\theta)$ on $Y$ for each $\nu.$
Denote by $V_{\pi,m}$ the subspace
of $V_{\pi}$ consisting of $P\cap K_m$--fixed
vectors.  Since $\pi$
is admissible, $V_{\pi,m}$ is finite dimensional.

The next three results are standard.  We include the proof of the first
two for completeness.  The third
 is a straightforward consequence of the Iwasawa decomposition.

\proclaim{Lemma 3.1}  $Y_m$ has a basis $\{f_j\}$
which satisfies the following properties:
\roster
\item If $\g\in\Gamma_m,$ then the non-zero vectors among
$\{f_j(\g)\}$ are a basis for $V_{\pi,m}.$
\item If $f_j$ is fixed, then $f_j(\g)\neq 0$
for some $\g\in\Gamma_m.$
\endroster
\endproclaim
\demo{Proof}
Suppose that $f\in Y_m.$
Then $f(pkk_1)=\pi(p)f(k),$ for all $p\in P\cap K,$
$k\in K,$ and $k_1\in K_m.$  Thus, $f$ is completely determined
by its values on $\Gamma_m.$  Fix $\g\in\Gamma_m,$
and let $p\in P\cap K_m.$
Since $\g^{-1} K_m\g= K_m,$
we have $\g^{-1}p\g\in K_m.$
Therefore,
$$f(\g)=f(\g\g^{-1}p\g)=f(p\g)=\pi(p)f(\g).$$
This says that $f(\g)$ is an element of $V_{\pi,m}.$ 
Fix a basis $\{v_{m,i}\}$ of $V_{\pi,m}.$
Let $f_{\g,i}:K\longrightarrow V_{\pi,m}$
be given by
$$f_{\g,i}(k)=\cases \pi(p)v_{m,i}&\text{ if }k=p\g k_1,\text{ for some }p\in P\cap K,k_1\in K_m,\\0&\text{otherwise.}\endcases$$
Then it is immediate that
$f_{\g,i}$ is a well-defined element of $Y_m.$  We claim
that $\{f_{\g,i}\}$ is a basis for $Y_m.$

Suppose $f\in Y_m.$
If $\g'\in\Gamma_m,\, p\in P\cap K,$ and $k\in K_m,$
then
$$f(p\g' k)=\pi(p)f(\g').$$  Since $f(\g')\in V_{\pi,m},$
$$f(\g')=\sum_{i}c_{\g',i}v_{m,i}.$$
This implies that
$$f(p\g' k)=\sum_i c_{\g',i}\pi(p)v_{m,i}=\sum_ic_{\g',i}f_{\g',i}(p\g' k).$$
Now, taking the collection $\{c_{\g,i}\}$ for all $\g\in\Gamma_m,$
and noting that $f_{\g,i}(p\g'k)=0$ for $\g\neq\g',$
$$f(p\g' k)=\sum_{\g,i} c_{\g,i}f_{\g,i}(p\g' k),$$
which says that $\{f_{\g,i}\}$ spans $Y_m.$

On the other hand, suppose that $\sum\limits_{\g,i}c_{\g,i}f_{\g,i}=0.$
Then, for any $\g'\in \Gamma_m,$
we have
$$\gather
\sum_{\g,i} c_{\g,i}f_{\g,i}(\g')=0,\\
\intertext{which implies that}
\sum_i c_{\g',i} f_{\g',i}(\g')=\sum_i c_{\g',i}v_{m,i}=0.
\endgather$$
But, since the $v_{m,i}$ are linearly independent, 
$c_{\g',i}=0,$ for each $\g'$ and $i.$  Thus, $\{f_{\g,i}\}$
are also linearly independent.
The collection $f_{\g,i}$ clearly has
properties (1) and (2).
\qed\enddemo

Denote by $X$
the space $Y\otimes\tilde V_\o.$  For each
$\nu\in\frak a^*_\C$ let $X_\nu=V(\nu,\pi,\theta)\otimes\tilde V_\o.$
For $f\in Y$ denote by $f_{\nu}$ the unique element of
$V(\nu,\pi,\T)$ satisfying $f_{\nu}|_K=f.$
Then $\{f_{\nu}\otimes\tilde v|f\in Y,\tilde v\in \tilde V_\o\}$
spans $X_\nu.$  
Recall that $\Pi_{\nu}$ can be realized on $Y$ via
$\Pi_{\nu}(g)f=\left[\Pi_{\nu}(g)f_{\nu}\right]|_K.$  
This gives the context
in which we discuss the holomorphicity of the
map $\nu\mapsto\Pi_{\nu}(g)f_{\nu}$ for a fixed choice of $g$
and $f.$

\proclaim{Lemma 3.2}  Fix  $g\in G,\, f\in Y$ and $\tilde v\in \tilde V_\o.$
Then the function
$\nu\mapsto\Pi_{\nu}(g)f\otimes\tilde v$ is a regular function
from  $\frak a_\C^*$ to $X.$
\endproclaim

\demo{Proof}
Choose $m_0$ so that $f\in Y_{m_0},$ and choose $m>m_0$
satisfying $g^{-1}K_mg\subset K_{m_0}.$ Then $f\in Y_m$
and, for all $\nu\in\frak a^*_{\C}$ and $k\in K_m,$
$$\Pi_{\nu}(k)(\Pi_{\nu}(g)f_{\nu})(x)=f_{\nu}(xgg^{-1}kg)=f_{\nu}(xg)=\Pi_{\nu}(g)f_{\nu}(x),$$
which says that $\Pi_{\nu}(g)f\in Y_m$ for all $\nu.$
Now, by Lemma 3.1,
$$\Pi_{\nu}(g)f=\sum\limits_{\g,i}c_{\g,i}(\nu) f_{\g,i}$$
for a unique choice of $c_{\g,i}(\nu)\in \C.$
It suffices to show that $c_{\g,i}:\frak a_{\C}^*\longrightarrow \C$ is
holomorphic.  
Fix $\g'\in\Gamma_m.$  Then $\g' g=p\g'' k,$
for some $p\in P,$ $\g''\in\Gamma_m,$ and $k\in K_m.$
Then
$$\align
\Pi_\nu(g)f(\g')&=q_F^{<\nu,H_P(p)>}\d_P^{1/2}(p)\pi(p)f(\g'' k)\\
&=q_F^{<\nu,H_P(p)>}\d_P^{1/2}(p)\pi(p)\sum\limits_{\g,i}c_{\g,i}(\nu)f_{\g,i}(\g'')\\
&=q_F^{<\nu,H_P(p)>}\d_P^{1/2}(p)\pi(p)\sum\limits_i c_{\g'',i}(\nu)v_i.
\endalign$$
Set $c'_{\g'',i}(\nu)=q_F^{<-\nu,H_p(p)>}\d_P^{-1/2}(p)c_{\g'',i}(\nu).$
Then
$$\pi(p)f(\g'')=\sum\limits_i c_{\g'',i}(\nu) v_i,$$
for all $\nu.$  Since the left hand side in the
equation above is independent of $\nu$
and the $v_i$ are linearly independent,
$c'_{\g'',i}(\nu)$ is constant for each $i.$  This implies that 
$c_{\g'',i}(\nu)$ is holomorphic.
\qed\enddemo

From now on we need to distinguish between a Weyl group
element $\tilde w\in W(\G,\A),$ for some torus $\A,$
and a representative $w\in N_{G}(A)$ for $\tilde w.$
Let $\tilde w_{\T}=\tilde w_{l,\D} \tilde w_{l,\T},$
where $\tilde w_{l,\D}$ is the longest element of the Weyl
group $W(\G,\TT),$ and $\tilde w_{l,\T}$ is the longest element of
$W(\G,\A_\T).$  Fix a representative $w_{\T}$
for $\tilde w_{\T}$ with $w_{\T}\in K.$
Note that $\tilde w_\T(\T)\subset \D.$
Now let $\M'=\M_{\tilde w_\T(\T)}=w_\T \M_\T w_{\T}^{-1}.$
Then $\M'$ is a standard Levi subgroup of $\G.$
Let $\N'$ be the standard unipotent subgroup of $\U$
so that $\PP'=\M'\N'$ is a standard parabolic subgroup of $\G.$
Since $\M'\simeq\M,$ we have $\U_\ell\supset\N'.$

\proclaim{Lemma 3.3}
For each $m>0,$ we have $w_\T^{-1}N'\cap Pw_\T^{-1}K_m$ is compact.
\endproclaim

For $m\in M_\t,$ let $\d(m)=(\d_P^{-1/2}\d_{P'})(m).$ Let
$X_{R_\t,\o,\nu,\theta}$
be the subspace spanned by functions of the form
$$\Pi_{\nu}(mu)f\otimes \tilde v-\d(m)\t(u)f\otimes\tilde\o(m^{-1})\tilde v,$$
for $m\in M_\t,$ $u\in U_\ell,$ $f\in Y,$ and
$\tilde v\in V_{\tilde\o}.$
Then 
a non-zero functional $\L$
on $X$ is a $(\d\o_\t)$\snug-Bessel functional for $\Pi_{\nu}$
if and only if $\L|_{X_{R_\t,\o,\nu,\theta}}\equiv 0.$  
By Theorem 2.1 the space of
such functionals is   one-dimensional.
Thus $X/X_{R_\t,\o,\nu,\theta}$ is one-dimensional.

The construction of this functional will
be obtained by taking a direct limit of functionals
given by integrating over compact subsets of $N'.$
We  show that such a limit exists and is not
identically zero.  Moreover, we show that there is
a function in $X$ which is a complement to 
$X_{R_\chi,\omega,\nu,\theta}$ for all $\nu$.
This will give the regularity condition
necessary to apply Bernstein's Theorem and
to obtain the holomorphicity of the functional.

Now let us fix a Whittaker functional for $\s$ and a Bessel functional
for $\o.$ 
(Actually, for notational convenience, we twist $\o$
by $\d_{R\t}^{-1/2}.$) That is, suppose that
$\l_{\t}:V_{\pi}\otimes \tilde V_\o\longrightarrow\C$
satisfies
$$\l_{\t}((\s(u_1)\otimes\tau(mu_2))(v_1\otimes v_2)\otimes\tilde v)=
\t_1(u_1)\t'(u_2)\l_{\t}(v_1\otimes v_2\otimes\tilde\o(m^{-1})\tilde v)$$
for all $u_1\in U_{\ell}\cap G_1,$ $u_2\in U_{\ell}\cap G(m),$
and $m\in M_\t.$
Let $\O$ be a compact subgroup of $N'.$
Define a functional on $X$ by
$$\l_{\pi,\nu,\T}^{\O}(f\otimes\tilde v)=\int\limits_{\O}\l_\t(\Pi_{\nu}(w_\T^{-1}u)f_\nu(e)\otimes\tilde v)\t(u)^{-1}\,du.\tag 3.1$$
This functional depends on the choice of the
representative $w_{\T}$ for $\tilde w_{\T}.$

Since $N'$ is exhausted by compact subgroups, the
compact subgroups of $N'$ form a directed set.
The following Lemma was suggested to the authors by Prof.\ Steve Rallis.

\proclaim{Lemma 3.4}
For every $f\otimes\tilde v\in X,$ the limit
$\lim\limits_{\O}\l_{\pi,\nu,\T}^{\O}(f\otimes\tilde v)$
exists, where the limit is the direct limit
taken over all compact subgroups of $N'.$
\endproclaim

\demo{Proof} Fix $f\otimes\tilde v\in X$.
Let $\Omega_r$ be the subset of $N'$ of elements with all
entries of absolute value at most $q_F^r$.  It suffices
to show that there is an $r$ sufficiently large 
such that if $\Omega\supset
\Omega_r$ then
$\l_{\pi,\nu,\T}^{\O}(f\otimes\tilde v)=
\l_{\pi,\nu,\T}^{\O_r}(f\otimes\tilde v).$ 
This follows since, for $r$ sufficiently
large, if $\gamma\in N'\setminus\O_r$ then
$$\int\limits_{\O_r\gamma}\l_\t(\Pi_{\nu}(w_\T^{-1}u)
f_\nu(e)\otimes\tilde v)\t(u)^{-1}\,du=0.$$
Indeed, there is a subgroup $K\subset\O_r$ such that
$f_\nu(w_\T^{-1}ku)=f_\nu(w_\T^{-1}u)$ 
for $k\in K$, $u\in \O_r\gamma$,   
but such that $\t$ is not identically 1 on $K$. 
For one may choose $K$ such that the relevant minors of $w_\T^{-1}u$ 
and $w_\T^{-1}ku$ are highly congruent, and use the Iwasawa 
decomposition.
\qed
\enddemo

Define a functional on $X$ by
$$\L_\t(\nu,\pi,\T)(f\otimes\tilde v)=\lim\limits_\O \l_{\pi,\nu,\T}^{\O}(f\otimes\tilde v).\tag 3.2$$
Again, this functional depends on the choice of
$w_{\T}.$

\proclaim{Proposition 3.5}
Let $\L_\t(\nu,\pi,\T)$ be defined as in {\rm (3.2)},
and extend $\L_\t$ to $X_\nu$ by the section
$f\otimes \tilde v\mapsto f_\nu\otimes \tilde v$.
Then  $\L_\t(\nu,\pi,\T)$ defines a non-zero
$\d\o_\chi$--Bessel functional for $\Pi_{\nu}.$\qed
\endproclaim
\demo{Proof} Suppose that $u_1\in U_\ell.$  Since $U_\ell\subset P',$
we can write $u_1=m_1n_1,$ with $m_1\in M'\cap U_\ell,$ and $n_1\in N'.$
Suppose first that $u_1=n_1\in N'.$  Since $N'$ is
exhausted by compact subgroups, we can choose $\O_0$ compact
with $n_1\in\O_0.$   If $\O_0\subset \O,$ then 
$$\gather
\l_{\pi,\nu,\T}^{\O}(\Pi_\nu(n_1)f\otimes\tilde v)=\int\limits_{\O} \l_\t(f_\nu(w_\T^{-1}un_1)\otimes\tilde v)\t^{-1}(u)\, du\\
=\int\limits_{\O}\l_{\t}(f_\nu(w_\T^{-1}u)\otimes\tilde v)\t^{-1}(un_1^{-1})\,du\\
=\t(n_1)\l_{\pi,\nu,\T}^{\O}(f\otimes\tilde v).
\endgather$$
Therefore, 
$$\L_\t(\nu,\pi,\T)(\Pi_\nu(n_1)f\otimes\tilde v)=\t(n_1)\L_\t(\nu,\pi,\T)(f\otimes\tilde v).$$
If $u=m_1\in U_\ell\cap M',$ then since $\t|_{G_1\cap U_\ell}$
is $\psi$--generic,
$\t^{w_\T}(m_1)=\t(m_1).$  Thus,
$$\gather
\l_{\pi,\nu,\T}^{\O}(\Pi_\nu(m_1)f\otimes\tilde v)=\int\limits_{\O}\l_\t(f_\nu(w_\T^{-1}um_1)\otimes\tilde v)\t^{-1}(u)\, du\\
=\int\limits_{\O}\l_\t(f_\nu(w_\T^{-1}m_1w_\T w_\T^{-1}m_1^{-1}um_1)\otimes\tilde v)\t^{-1}(u)\, du\\
=\int\limits_{\O}\l_\t(\pi(w_\T^{-1}m_1 w_\T)f_\nu(w_\T^{-1}m_1^{-1}um_1)\otimes\tilde v)\t^{-1}(u)\, du\\
=\t(m_1)\int\limits_{m_1^{-1}\O m_1}\l_{\t}(f_\nu(w_\T^{-1}u)\otimes\tilde v)\t^{-1}(u)\, du\\
=\t(m_1)\l_{\pi,\nu,\T}^{m_1^{-1}\O m_1}(f\otimes\tilde v).
\endgather$$
Therefore, 
$$\L_\t(\nu,\pi,\T)(\Pi_\nu(m_1)f\otimes\tilde v)=\t(m_1)\L_\t(\nu,\pi,\T)(f\otimes\tilde v).$$

Similarly, if $m\in M_\t\subset M',$
then
$$\gather
\l_{\pi,\nu,\T}^{\O}(\Pi_\nu(m)f\otimes\tilde v)=\int\limits_{\O}\l_{\t}(f_\nu(w_\T^{-1} u m)\otimes\tilde v)\t^{-1}(u)\, du\\
=\int\limits_{\O}\l_\t(\pi(w_\T^{-1}mw_{\T})\d_{P}^{1/2}(w_\T^{-1}mw_\T)
f_\nu(w_\T^{-1}m^{-1}um)\otimes\tilde v)\t^{-1}(u)\, du\\
=\d_P^{1/2}(w_\T^{-1}mw_\T)\int\limits_{m^{-1}\O m} \l_\t(f_\nu(w_\T^{-1}m^{-1}um)\otimes\tilde\o(m^{-1})\tilde v)\t^{-1}(u)\, du\\
=\d_P^{-1/2}\d_{P'}(m)\int\limits_{m^{-1}\O m}\l_{\t}(f_\nu(w_\T^{-1} u)\otimes\tilde\o(m^{-1})\tilde v)\t^{-1}(mum^{-1})\, du\\
=\d(m)\l_{\pi,\nu,\T}^{m^{-1}\O m}(f\otimes\tilde\o(m^{-1})\tilde v).
\endgather$$
Taking the limit on $\O$ on the right and left sides of
the above equation completes
the proof that
$\L_\t(\nu,\pi,\T)$ is a Bessel functional for
$\Pi_\nu$ with respect to the representation $\d(m)\o_\t.$

It remains to show that $\L_\t(\nu,\pi,\T)$ is not identically zero.
Let $\bar{\PP'}$ be the parabolic opposite to $\PP'$.  
Then $\bar{\PP'}=w_\T\PP w_\T^{-1}.$
By Lemma 3.4,
$\bar P' K_m$ is compact, and if $pw_\T^{-1}k\in Pw_\T^{-1}K_m\cap N',$
then in fact $p\in P\cap K_m.$
Choose a $v\in V_{\pi}$ and $\tilde v\in \tilde V_\o$
such that $\l_\t(v\otimes\tilde v)\neq 0.$
Choose $m>>0$ such that
$v\in V_{\pi,m}$ and such that  $\t|{N'\cap\bar P' K_m}\equiv 1.$
Consider the function in $Y$ defined by
$$f_0(k)=\cases \pi(p)v&\text{if }k=pw_\T^{-1}k_1, p\in P\cap K,k_1\in K_m\\
0&\text{otherwise.}\endcases\tag 3.3
$$
Then
$$\gather
\L_\t(\nu,\pi,\T)(f_0\otimes\tilde v)=\int\limits_{N'\cap \bar P' K_m}\l_\t(f_\nu(w_\T^{-1}u)\otimes\tilde v)\, du\\
=\l_\t(v\otimes\tilde v)\, |N'\cap\bar P'K_m|\neq 0.
\endgather$$
Thus, $\L_\t(\nu,\pi,\T)$ is non-zero, and $f_0$
is a complement to $X_{R_\t,\o,\nu,\T}$ for all $\nu.\qed$
\enddemo

Suppose $r=mu\in R_{\t},$ $f\in Y,$
and $\tilde v\in \tilde V_{\o}.$
Define an $X$--valued
function on $\frak a_{\C}^*$ by
$$x_{r,f,\tilde v,\theta}(\nu)=\Pi_{\nu}(r)(f)\otimes\tilde v-\d(m)\t(u)(f\otimes\tilde\o(m^{-1})\tilde v).$$

\proclaim{Theorem 3.6}
The function $\nu\mapsto \L_\t(\nu,\pi,\T)(x)$ is holomorphic for each
$x\in X.$
\endproclaim

\demo{Proof}  We will apply Banks's extension of Bernstein's Theorem.
Let 
$$\Cal R=\{(r,f\otimes\tilde v)|r\in R_{\t}, f\in Y,\tilde v\in V_{\tilde\o}\}\cup\{*\}.$$
For $\a=(r,f\otimes\tilde v)\in\Cal R,$ we let $x_\a(\nu)=x_{r,f,\tilde v,\T}(\nu)$
in $X$ and let $c_\a(\nu)=0.$  Fix $m,$ $\tilde v,$ $v,$
and $f_0$ as in (3.3). 
For $\a=*,$ we set $x_*(\nu)=f_0\otimes\tilde v$
and $c_*(\nu)=|N'\cap\bar P' K_m\,|\,\l_\t(v\otimes \tilde v).$
Now for every $\nu\in\frak a_{\C}^*,$ we consider the systems of
equations in $X\times\C$ given by:
$$\Xi(\nu)=\left\{(x_\a(\nu),c_\a(\nu))|\,\a\in \Cal R\right\}.$$
By Lemma 3.2, the function $\nu\mapsto x_\a(\nu)$ is 
holomorphic for each $\a$ of the form $(r,f\otimes\tilde v).$
For $\a=*,$ the function
$x_\a(\nu)=f_0\otimes\tilde v$ is constant on $\frak a^*_\C.$
Note that each $c_\a$ is constant, hence holomorphic as well.

Now,  for each $\nu$ the functional $\L_\t(\nu,\pi,\T)$ is a solution
to the system $\Xi(\nu).$  Moreover, such a solution
is unique by the results of Section 2.  Thus, Banks's
extension  of Bernstein's theorem\Lspace \Lcitemark Ban\Rcitemark \Rspace{} implies that $\nu\mapsto\L(\nu,\pi,\T)(f\otimes\tilde v)$
is holomorphic for all choices of $f$ and $\tilde v.$\qed\enddemo

We turn to the question of local coefficients.
Let $\tilde w\in W,$ and fix
a representative $w$ for $\tilde w$ with $w\in K.$
We recall that the intertwining operator
$$A(\nu,\pi,w):V(\nu,\pi,\T)\rightarrow V(\tilde w(\nu),\tilde w\pi,\tilde w(\T))$$
is defined for $\nu >>0$ by
$$A(\nu,\pi,w)f(g)=\int\limits_{N_{\tilde w}} f(w^{-1}ng)\, dn,\tag 3.4$$
where $\N_{\tilde w}=\U\cap w\bar\N w^{-1},$ and
$\bar{\N}$ is the unipotent radical opposite to $\N.$
Then $A(\nu,\pi,w)$ is defined on all of $\frak a_{\C}^*$ by
analytic continuation.  Note that the intertwining operator
depends on the choice of $w$ representing $\tilde w.$

We also recall the Langlands decomposition of the intertwining
operator, described  in Lemma 2.1.2 of\Lspace \Lcitemark Shaa\Rcitemark \Rspace{}.
For the convenience of the reader,
let us restate this here.  
For two associate subsets $\theta$ and $\theta'$ of $\Delta,$
we let
$$W(\theta,\theta')=\left\{\tilde w\in W\, |\, \tilde w\theta=\theta'\right\}.$$
\proclaim{Lemma 3.7 (Langlands (see\Lspace \Lcitemark Shaa\LIcitemark{}, Lemma 2.1.2\RIcitemark \Rcitemark \Rspace{}))}
Suppose that $\T,\T'\subset \D$ are associate.  Let $\tilde w\in W(\T,\T').$
Then there exists a family $\T_1,\T_2,\dots,\T_n\subset\D$
so that
\roster
\item $\T_1=\T$ and $\T_n=\T';$
\item For each $1\leq i\leq n$ there is a root
$\a_i\in\D\setminus\T_i$ so that $\T_{i+1}$
is the conjugate of $\T_i$ in $\D_i=\T_i\cup\{\a_i\};$
\item For each $1\leq i\leq n-1,$ we let $\tilde w_i=\tilde w_{\ell,\D_i}\tilde w_{\ell,\T_i}$
in $W(\T_i,\T_{i+1}).$  Then $\tilde w=\tilde w_{n-1}\dots \tilde w_1;$
\item Set $\tilde w_1'=\tilde w,$ and $\tilde w'_{i+1}= \tilde w_i'\tilde w_i^{-1}$
for $1\leq i\leq n-1.$  Then
$\tilde w_n'=1$ and
$$\frak n_{\tilde w_i'}=\frak n_{\tilde w_i}\oplus\text{{\rm Ad}}(w_i^{-1})\frak n_{\tilde w'_{i+1}}.$$
Here $\frak n$ is the Lie algebra of $\N.$
\endroster
\endproclaim

Let $\T_{*}\subset\T$ and let $\rho$ be
an irreducible supercuspidal representation 
of $M_{\T_*}$. 
If $\rho$ is generic, then Rodier's Theorem implies that 
there is a unique constituent $\pi$ of
$\Ind_{P_{\T_*}}^{M}(\rho)$ which is generic 
with compatible character.
For this constituent, Shahidi proved that there is a 
complex number
$C_\t(\nu,\pi,\T,w)$ which satisfies
$$\L_\t(\nu,\pi,\T)=C_\t(\nu,\pi,\T,w)\L_\t(\tilde w\nu,\tilde w\pi,\tilde w\T)A(\nu,\pi,w),$$
where $\L_\t$ is the Whittaker functional.
Moreover, the function $\nu\mapsto C_\t(\nu,\pi,\T,w)$
is a meromorphic function on $(\frak a_\T)^*_\C.$
The value of the local coefficient depends on the choice
of representative $w$ for $\tilde w.$

Now suppose that $\rho$ is any irreducible supercuspidal
which has a minimal Bessel model of a particular type.
Then we prove a similar result for the constituent
$\pi$ of
$\Ind_{P_{\T_*}}^{M}(\rho)$ 
which has a Bessel model of 
compatible type; such a
constituent exists and is unique by
Theorem 2.1.
%Below we prove a similar result with one additional assumption,
%which we now describe.
%Write $\rho=\rho_1\otimes\rho_2,$ where $\rho_1$
%is an irreducible supercuspidal representation of
%$$GL_{a_1}(\eurm k)\times\dots\times GL_{a_s}(\eurm k),$$
%and $\rho_2$ is an irreducible supercuspidal representation of $G(m')$
%for some $m'\leq m.$  Then $\rho_1$ is generic and $\Cal B(\rho_2)\leq \Cal B(\tau).$
%Suppose that $\Cal B(\rho_2)=\Cal B(\tau).$
%Note that $I(\nu,\pi,\T)\hookrightarrow I(\tilde\nu,\rho,\T_*),$
%when $\tilde\nu$ restricts to $\nu.$  Furthermore, Theorem 2.1  
%guarantees that both of these representations have a unique
%Bessel model of a given type of rank $\Cal B(\rho_2).$  Thus, the
%constituent of $I(\nu,\pi,\T)$ which carries a minimal Bessel model,
%is also the unique such constituent of $I(\tilde\nu,\rho,\T_*).$
%It is in this case that we can attach a local coefficient to
%$\pi.$

\proclaim{Theorem 3.8} Let $\T$ and $\T'$ be associate subsets of $\D.$
Let $\T_{*}\subset\T$ and let $\rho$ be
an irreducible supercuspidal representation  
of $M_{\T_*}$.  Suppose that $\rho$ has an
$\omega_\chi$--Bessel model which is minimal. 
Let $\pi$ be the constituent of 
$\Ind_{P_{\T_*}}^{M}(\rho)$ such that 
$\pi$ has an $\omega_\chi^{w_0}$-Bessel model,
as in Theorem 2.1.
For each $\tilde w\in W(\T,\T')$ fix a representative $w$
for $\tilde w.$  Then  there is a complex number $C_\t(\nu,\pi,\T,w)$
so that
$$\L_\t(\nu,\pi,\T)=C_\t(\nu,\pi,\T,w)\L_\t(\tilde w\nu,\tilde w\pi,\tilde w\T)A(\nu,\pi,w).\tag 3.5$$
Moreover, the function $\nu\mapsto C_\t(\nu,\pi,\T,w)$ is meromorphic
on $\frak a_\C^*,$ and depends only on the class of $\pi$
and the choice of $w.$
\endproclaim

\demo{Proof}  We first show how to define $C_\t(\tilde\nu,\pi,\T_*,w)$
for $\tilde\nu\in (\frak a_{\T_*})^*_{\C}.$
By\Lspace \Lcitemark Sil\LIcitemark{}, Theorem 5.4.3.7\RIcitemark \Rcitemark \Rspace{} the representation
$I(\tilde\nu,\rho,\T_*)$ is irreducible unless the Plancherel measure 
$\mu(\tilde\nu,\rho)=0$ and $(\tilde\nu,\rho)$ is fixed by a nontrivial
element of the Weyl group $W_{\T_*}$ (i.e., is singular).
Thus, on an open dense subset of $(\frak a_{\T_*})^*_\C$
the representation $I(\tilde\nu,\rho,\T_*)$ is irreducible, and so $\L_\t(\tilde w\tilde\nu,\tilde w\rho,\tilde w\T_*)A(\tilde\nu,\rho,w)$
defines a non-zero  Bessel functional on 
$V(\tilde\nu,\rho,\T_*)\otimes\tilde V_\o.$
By the uniqueness of such a functional 
(Theorem/Conjecture 1.4), we get the
existence of $C(\tilde\nu,\rho,\T_*,w)$ 
satisfying
$$\L_\t(\tilde\nu,\rho,\T_*)=C_\t(\tilde\nu,\rho,\T_*,w)\L_\t(\tilde w\tilde\nu,\tilde w\rho,\tilde w\T_*)A(\tilde\nu,\rho,w)$$
on the open dense subset.
Moreover, it is holomorphic there since both 
$$\L_\t(\tilde w\tilde\nu,\tilde w\rho,\tilde w\T_*)$$
and $A(\tilde\nu,\rho,w)$ are holomorphic there.  Thus, $C(\tilde\nu,\rho,\T_*,w)$
extends to a meromorphic function on $(\frak a_{\T_*})^*_\C.$
Now, write $\tilde w=\tilde w_{n-1}\dots \tilde w_1$ as in Lemma  3.7.
Since $C_\t(\tilde\nu,\rho,\T_*,w)$ is now defined,
it admits a factorization compatible with the decomposition
of the intertwining operators given in Lemma 3.7.
(See corollary 3.9.)
This implies that on an open dense subset of $\nu\in(\frak a_\T)^*_\C,$
the local coefficient $C_\t(\nu,\rho,\T_*,w)$ may be defined
by the equation
$C_\t(\nu,\rho,\T_*,w)=C_\t(\tilde\nu,\rho,\T_*,w),$
where $\tilde\nu$ is the restriction of $\nu$ to $(\frak a_{\T_*})_\C.$
Suppose that, for some $\nu$ in this open dense subset,
$\L_\t(\tilde w\nu,\tilde w\pi,\tilde w\T)A(\nu,\pi,w)$ was the
zero functional.  Then, by inducing in stages and using the
discussion preceding this Theorem,
we would conclude that
$\L_\t(\tilde w\tilde\nu,\tilde w\rho,\tilde w\T_*)A(\tilde\nu,\rho,w)$
is also zero.  However, since $C_\t(\tilde\nu,\rho,\T_*,w)$
is defined there, this would be  a contradiction. 
Thus we may define $C(\nu,\pi,\T,w)$ by the relation (3.5) on this
open dense subset, and we have $C(\nu,\pi,\T,w)=C(\tilde\nu,\rho,\T_*,w).$
Since $A(\nu,\pi,w)$ has a meromorphic continuation
to $(\frak a_\T)^*_\C,$ and $\L_\t(\tilde w\nu,\tilde w\pi,\tilde w\T)$
is holomorphic on $(\frak a_\T)^*_\C,$
the function $\nu\mapsto C_\t(\nu,\pi,\T,w)$ must
have a meromorphic continuation.
\qed
\enddemo

\proclaim{Corollary 3.9}  Let the notation be as in Lemma 3.7
and Theorem 3.8.
Let  $\pi_1=\pi,$ and $\nu_1=\nu.$  For each $i,$ $2\leq i\leq n-1,$
set $\pi_i=\tilde w_i\pi_{i-1},\nu_i=\tilde w\nu_{i-1}.$  Then the local
coefficient factors as
$$C_\t(\nu,\pi,\T,w)=\prod_{i=1}^{n-1} C_\t(\pi_i,\T_i,w_i).$$
\endproclaim
\demo{Proof}  Let $f_1=f\in V(\nu,\pi,\T)$  and for $2\leq i\leq n-1,$
let $f_i=A(\nu_{i-1},\pi_{i-1},w_{i-1})f_{i-1}.$  Then
$$\align
\L_\t(\nu_i,\pi_i,\T_i)f_i=&C_\t(\nu_i,\pi_i,\T_i,w_i)\L_\t(\nu_{i+1},\pi_{i+1},\T_{i+1})\\
&\cdot A(\nu_i,\pi_i,w_i)f_i,
\endalign$$
for each $1\leq i\leq n-1.$
The corollary now follows immediately from Lemma 3.7
and iteration of the above equality.\qed
\enddemo

\vfil\eject
\Refs

\message{REFERENCE LIST}

\bgroup\Resetstrings%
\def\Ecnt{0}\def\acnt{0}%
\def\Ftest{ }\def\Fstr{Ban}%
\def\Atest{ }\def\Astr{W. B\bgroup\Smallcapsfont ANKS\egroup{}}%
\def\Ttest{ }\def\Tstr{Exceptional representations on the metaplectic group}%
\def\Otest{ }\def\Ostr{preprint}%
\Refformat\egroup%

\bgroup\Resetstrings%
\def\Ecnt{0}\def\acnt{0}%
\def\Ftest{ }\def\Fstr{BeZa}%
\def\Atest{ }\def\Astr{I.N. B\bgroup\Smallcapsfont ERNSTEIN\egroup{}%
  \Aand A.V. Z\bgroup\Smallcapsfont ELEVINSKY\egroup{}}%
\def\Ttest{ }\def\Tstr{Representations of the group $GL(n,F)$ where $F$ is a local non-archimedean field}%
\def\Jtest{ }\def\Jstr{Russian Math. Surveys}%
\def\Vtest{ }\def\Vstr{33}%
\def\Ptest{ }\def\Pstr{1--68}%
\def\Dtest{ }\def\Dstr{1976}%
\Refformat\egroup%

\bgroup\Resetstrings%
\def\Ecnt{0}\def\acnt{0}%
\def\Ftest{ }\def\Fstr{BeZb}%
\def\Atest{ }\def\Astr{I.N. B\bgroup\Smallcapsfont ERNSTEIN\egroup{}%
  \Aand A.V. Z\bgroup\Smallcapsfont ELEVINSKY\egroup{}}%
\def\Ttest{ }\def\Tstr{Induced representations of reductive $p$--adic groups. I}%
\def\Jtest{ }\def\Jstr{Ann. Sci. \'Ecole Norm. Sup. (4)}%
\def\Vtest{ }\def\Vstr{10}%
\def\Ptest{ }\def\Pstr{441--472}%
\def\Dtest{ }\def\Dstr{1977}%
\def\Astr{\Underlinemark}%
\Refformat\egroup%

\bgroup\Resetstrings%
\def\Ecnt{0}\def\acnt{0}%
\def\Ftest{ }\def\Fstr{Ber}%
\def\Atest{ }\def\Astr{J. B\bgroup\Smallcapsfont ERNSTEIN\egroup{}}%
\def\Ttest{ }\def\Tstr{Letter to I.I. Piatetski-Shapiro, Fall 1985}%
\def\Otest{ }\def\Ostr{to appear in a book by J. Cogdell and I.I. Piatetski-Shapiro}%
\Refformat\egroup%

\bgroup\Resetstrings%
\def\Ecnt{0}\def\acnt{0}%
\def\Ftest{ }\def\Fstr{Bru}%
\def\Atest{ }\def\Astr{F. B\bgroup\Smallcapsfont RUHAT\egroup{}}%
\def\Ttest{ }\def\Tstr{Sur les repr\'esentations induites des groupes de Lie}%
\def\Jtest{ }\def\Jstr{Bull. Soc. Math. France}%
\def\Dtest{ }\def\Dstr{1956}%
\def\Vtest{ }\def\Vstr{84}%
\def\Ptest{ }\def\Pstr{97--205}%
\Refformat\egroup%

\bgroup\Resetstrings%
\def\Ecnt{0}\def\acnt{0}%
\def\Ftest{ }\def\Fstr{Cas}%
\def\Atest{ }\def\Astr{W. C\bgroup\Smallcapsfont ASSELMAN\egroup{}}%
\def\Ttest{ }\def\Tstr{Introduction to the theory of admissible representations of $p$--adic reductive groups}%
\def\Otest{ }\def\Ostr{preprint}%
\Refformat\egroup%

\bgroup\Resetstrings%
\def\Ecnt{0}\def\acnt{0}%
\def\Ftest{ }\def\Fstr{GPR}%
\def\Atest{ }\def\Astr{D. G\bgroup\Smallcapsfont INZBURG\egroup{}%
  \Acomma I. P\bgroup\Smallcapsfont IATETSKII\egroup{}-S\bgroup\Smallcapsfont HAPIRO\egroup{}%
  \Aandd S. R\bgroup\Smallcapsfont ALLIS\egroup{}}%
\def\Ttest{ }\def\Tstr{$L$--functions for the orthogonal group}%
\def\Otest{ }\def\Ostr{preprint}%
\Refformat\egroup%

\bgroup\Resetstrings%
\def\Ecnt{0}\def\acnt{0}%
\def\Ftest{ }\def\Fstr{Gola}%
\def\Atest{ }\def\Astr{D. G\bgroup\Smallcapsfont OLDBERG\egroup{}}%
\def\Ttest{ }\def\Tstr{Reducibility of induced representations for $Sp(2n)$ and $SO(n)$}%
\def\Jtest{ }\def\Jstr{Amer. J. Math.}%
\def\Vtest{ }\def\Vstr{116}%
\def\Dtest{ }\def\Dstr{1994}%
\def\Ptest{ }\def\Pstr{1101--1151}%
\Refformat\egroup%

\bgroup\Resetstrings%
\def\Ecnt{0}\def\acnt{0}%
\def\Ftest{ }\def\Fstr{Golb}%
\def\Atest{ }\def\Astr{D. G\bgroup\Smallcapsfont OLDBERG\egroup{}}%
\def\Ttest{ }\def\Tstr{$R$--groups and elliptic representations for unitary groups}%
\def\Jtest{ }\def\Jstr{Proc. Amer. Math. Soc.}%
\def\Vtest{ }\def\Vstr{123}%
\def\Dtest{ }\def\Dstr{1995}%
\def\Ptest{ }\def\Pstr{1267--1276}%
\def\Astr{\Underlinemark}%
\Refformat\egroup%

\bgroup\Resetstrings%
\def\Ecnt{0}\def\acnt{0}%
\def\Ftest{ }\def\Fstr{Har}%
\def\Atest{ }\def\Astr{H\bgroup\Smallcapsfont ARISH\egroup{}-C\bgroup\Smallcapsfont HANDRA\egroup{}}%
\def\Ttest{ }\def\Tstr{Harmonic analysis on reductive $p$--adic groups}%
\def\Dtest{ }\def\Dstr{1973}%
\def\Jtest{ }\def\Jstr{Proc. Sympos. Pure Math.}%
\def\Itest{ }\def\Istr{AMS}%
\def\Ctest{ }\def\Cstr{Providence, RI}%
\def\Vtest{ }\def\Vstr{26}%
\def\Ptest{ }\def\Pstr{167--192}%
\Refformat\egroup%

\bgroup\Resetstrings%
\def\Ecnt{0}\def\acnt{0}%
\def\Ftest{ }\def\Fstr{HeR}%
\def\Atest{ }\def\Astr{E. H\bgroup\Smallcapsfont EWITT\egroup{}%
  \Aand K. R\bgroup\Smallcapsfont OSS\egroup{}}%
\def\Ttest{ }\def\Tstr{Abstract Harmonic Analysis I}%
\def\Itest{ }\def\Istr{Springer-Verlag}%
\def\Ctest{ }\def\Cstr{New York--Heidelberg--Berlin}%
\def\Dtest{ }\def\Dstr{1979}%
\Refformat\egroup%

\bgroup\Resetstrings%
\def\Ecnt{0}\def\acnt{0}%
\def\Ftest{ }\def\Fstr{Lana}%
\def\Atest{ }\def\Astr{R.P. L\bgroup\Smallcapsfont ANGLANDS\egroup{}}%
\def\Ttest{ }\def\Tstr{Euler Products}%
\def\Otest{ }\def\Ostr{Yale University}%
\def\Dtest{ }\def\Dstr{1971}%
\Refformat\egroup%

\bgroup\Resetstrings%
\def\Ecnt{0}\def\acnt{0}%
\def\Ftest{ }\def\Fstr{Lanb}%
\def\Atest{ }\def\Astr{R.P. L\bgroup\Smallcapsfont ANGLANDS\egroup{}}%
\def\Ttest{ }\def\Tstr{On the functional equation satisfied by Eisenstein series}%
\def\Itest{ }\def\Istr{Springer-Verlag}%
\def\Ctest{ }\def\Cstr{New York--Heidelberg--Berlin}%
\def\Stest{ }\def\Sstr{Lecture Notes in Mathematics }%
\def\Ntest{ }\def\Nstr{544}%
\def\Dtest{ }\def\Dstr{1976}%
\def\Astr{\Underlinemark}%
\Refformat\egroup%

\bgroup\Resetstrings%
\def\Ecnt{0}\def\acnt{0}%
\def\Ftest{ }\def\Fstr{Nov}%
\def\Atest{ }\def\Astr{M. N\bgroup\Smallcapsfont OVODVORSKY\egroup{}}%
\def\Ttest{ }\def\Tstr{New unique models of representations of unitary groups}%
\def\Jtest{ }\def\Jstr{Compositio Math.}%
\def\Vtest{ }\def\Vstr{33}%
\def\Dtest{ }\def\Dstr{1976}%
\def\Ptest{ }\def\Pstr{289--295}%
\Refformat\egroup%

\bgroup\Resetstrings%
\def\Ecnt{0}\def\acnt{0}%
\def\Ftest{ }\def\Fstr{Ral}%
\def\Atest{ }\def\Astr{S. R\bgroup\Smallcapsfont ALLIS\egroup{}}%
\def\Ttest{ }\def\Tstr{On certain Gelfand Graev models which are Gelfand pairs}%
\def\Otest{ }\def\Ostr{preprint}%
\Refformat\egroup%

\bgroup\Resetstrings%
\def\Ecnt{0}\def\acnt{0}%
\def\Ftest{ }\def\Fstr{Roda}%
\def\Atest{ }\def\Astr{F. R\bgroup\Smallcapsfont ODIER\egroup{}}%
\def\Ttest{ }\def\Tstr{Mod\`ele de Whittaker et caract\`eres de repr\'esentations}%
\def\Btest{ }\def\Bstr{Non Commutative Harmonic Analysis}%
\def\Itest{ }\def\Istr{Springer-Verlag}%
\def\Ctest{ }\def\Cstr{New York--Heidelberg--Berlin}%
\def\Stest{ }\def\Sstr{Lecture Notes in Math.}%
\def\Ntest{ }\def\Nstr{466}%
\def\Ptest{ }\def\Pstr{151--171}%
\Refformat\egroup%

\bgroup\Resetstrings%
\def\Ecnt{0}\def\acnt{0}%
\def\Ftest{ }\def\Fstr{Rodb}%
\def\Atest{ }\def\Astr{F. R\bgroup\Smallcapsfont ODIER\egroup{}}%
\def\Ttest{ }\def\Tstr{Whittaker models for admissible representations of reductive $p$--adic split groups}%
\def\Jtest{ }\def\Jstr{Proc. Sympos. Pure Math.}%
\def\Itest{ }\def\Istr{AMS}%
\def\Ctest{ }\def\Cstr{Providence, RI}%
\def\Vtest{ }\def\Vstr{26}%
\def\Ptest{ }\def\Pstr{425--430}%
\def\Dtest{ }\def\Dstr{1973}%
\def\Astr{\Underlinemark}%
\Refformat\egroup%

\bgroup\Resetstrings%
\def\Ecnt{0}\def\acnt{0}%
\def\Ftest{ }\def\Fstr{Shaa}%
\def\Atest{ }\def\Astr{F. S\bgroup\Smallcapsfont HAHIDI\egroup{}}%
\def\Ttest{ }\def\Tstr{On certain $L$--functions}%
\def\Jtest{ }\def\Jstr{Amer. J. Math.}%
\def\Vtest{ }\def\Vstr{103}%
\def\Ntest{ }\def\Nstr{2}%
\def\Ptest{ }\def\Pstr{297--355}%
\def\Dtest{ }\def\Dstr{1981}%
\Refformat\egroup%

\bgroup\Resetstrings%
\def\Ecnt{0}\def\acnt{0}%
\def\Ftest{ }\def\Fstr{Shab}%
\def\Atest{ }\def\Astr{F. S\bgroup\Smallcapsfont HAHIDI\egroup{}}%
\def\Ttest{ }\def\Tstr{Local coefficients as Artin factors for real groups}%
\def\Jtest{ }\def\Jstr{Duke Math. J.}%
\def\Vtest{ }\def\Vstr{52}%
\def\Dtest{ }\def\Dstr{1985}%
\def\Ptest{ }\def\Pstr{973--1007}%
\def\Astr{\Underlinemark}%
\Refformat\egroup%

\bgroup\Resetstrings%
\def\Ecnt{0}\def\acnt{0}%
\def\Ftest{ }\def\Fstr{Shac}%
\def\Atest{ }\def\Astr{F. S\bgroup\Smallcapsfont HAHIDI\egroup{}}%
\def\Ttest{ }\def\Tstr{On the Ramanujan conjecture and finiteness of poles for certain $L$--functions}%
\def\Jtest{ }\def\Jstr{Ann. of Math. (2)}%
\def\Vtest{ }\def\Vstr{127}%
\def\Dtest{ }\def\Dstr{1988}%
\def\Ptest{ }\def\Pstr{547--584}%
\def\Astr{\Underlinemark}%
\Refformat\egroup%

\bgroup\Resetstrings%
\def\Ecnt{0}\def\acnt{0}%
\def\Ftest{ }\def\Fstr{Shad}%
\def\Atest{ }\def\Astr{F. S\bgroup\Smallcapsfont HAHIDI\egroup{}}%
\def\Ttest{ }\def\Tstr{A proof of Langlands conjecture for Plancherel measures; complementary series for $p$--adic groups}%
\def\Jtest{ }\def\Jstr{Ann. of Math. (2)}%
\def\Dtest{ }\def\Dstr{1990}%
\def\Vtest{ }\def\Vstr{132}%
\def\Ptest{ }\def\Pstr{273--330}%
\def\Astr{\Underlinemark}%
\Refformat\egroup%

\bgroup\Resetstrings%
\def\Ecnt{0}\def\acnt{0}%
\def\Ftest{ }\def\Fstr{Sil}%
\def\Atest{ }\def\Astr{A.J. S\bgroup\Smallcapsfont ILBERGER\egroup{}}%
\def\Ttest{ }\def\Tstr{Introduction to Harmonic Analysis on Reductive $p$--adic Groups}%
\def\Itest{ }\def\Istr{Princeton University Press}%
\def\Ctest{ }\def\Cstr{Princeton, NJ}%
\def\Stest{ }\def\Sstr{Mathematical Notes}%
\def\Ntest{ }\def\Nstr{23}%
\def\Dtest{ }\def\Dstr{1979}%
\Refformat\egroup%

\endRefs
\enddocument
\end